\documentclass[12pt,leqno]{article}
\usepackage{amsmath}
\usepackage{amssymb}
\usepackage[matrix,arrow]{xy}

\numberwithin{equation}{section}

\CompileMatrices

\newcounter{rocount}
\newenvironment{rolist}{\begin{list}{(\roman{rocount})}{\usecounter{rocount}\setlength{\itemsep}{0.0cm}\setlength{\parsep}{0.0cm}}}{\end{list}}

\newenvironment{Rolist}{\begin{list}{(\Roman{rocount})}{\usecounter{rocount}\setlength{\itemsep}{0.0cm}\setlength{\parsep}{0.0cm}}}{\end{list}}

\newcounter{gracount}

\newtheorem{theor+}{Theorem}[subsection]
\newtheorem{propo+}[theor+]{Proposition}
\newtheorem{corol+}[theor+]{Corollary}
\newtheorem{lemma+}[theor+]{Lemma}
\newtheorem{obser+}[theor+]{Observation}
\newtheorem{defin+}[theor+]{Definition}
\newtheorem{examp+}[theor+]{Example}
\newtheorem{remar+}[theor+]{Remark}
\newtheorem{conje+}[theor+]{Conjecture}
\newtheorem{conve+}[theor+]{Convention}
\newtheorem{quest+}[theor+]{Question}
\newtheorem{pictu+}[theor+]{Picture}
\newtheorem{notes+}[theor+]{Notes}
\newtheorem{exerc+}[theor+]{Exercise}
\newtheorem{backg+}[theor+]{Background}

\newenvironment{theor}[1]{\begin{theor+}\label{#1}\slshape}{\end{theor+}}
\newenvironment{propo}[1]{\begin{propo+}\label{#1}\slshape}{\end{propo+}}
\newenvironment{corol}[1]{\begin{corol+}\label{#1}\slshape}{\end{corol+}}
\newenvironment{lemma}[1]{\begin{lemma+}\label{#1}\slshape}{\end{lemma+}}

\newenvironment{defin}[1]{\begin{defin+}\label{#1}\upshape}{\end{defin+}}
\newenvironment{examp}[1]{\begin{examp+}\label{#1}\slshape}{\end{examp+}}
\newenvironment{remar}[1]{\begin{remar+}\label{#1}\upshape}{\end{remar+}}

\newenvironment{notes}{\begin{notes+}\upshape}{\end{notes+}}

\newenvironment{demo}{{\bf\sl Proof:}}{\hfill$\blacksquare$}

\newcommand{\genabsrep}{\gamma}
\newcommand{\Zgenabsrep}{\Gamma}
\newcommand{\Egen}{E_\genabsrep}

\newcommand{\hoo}{\mbox{$*$}-homomorphism}
\newcommand{\hoos}{\mbox{$*$}-homomorphisms}
\newcommand{\Ad}{\operatorname{Ad}}
\newcommand{\Aut}{\operatorname{Aut}}
\newcommand{\CALrho}{\dot{\MZtau}}
\newcommand{\CCC}{{\mathbb C}}

\newcommand{\COM}{^c}

\newcommand{\CPAM}{completely positive contractive map}
\newcommand{\CPAMs}{completely positive contractive maps}
\newcommand{\Dilim}{{\displaystyle\lim_{\longrightarrow}}}
\newcommand{\Drho}{D_\MZtau}

\newcommand{\ETA}{\underline{\pmb{\eta}}}
\newcommand{\Ext}{\operatorname{Ext}}
\newcommand{\Hdo}{H_{3}}
\newcommand{\Htwo}{H_{2}}
\newcommand{\Hha}{H_{1}}
\newcommand{\Hom}{\operatorname{Hom}}
\newcommand{\Hquad}{H_{6}}
\newcommand{\Hwh}{H_{1}}
\newcommand{\Inn}{\operatorname{Inn}}

\newcommand{\KKK}{{\mathcal K}}
\newcommand{\KKhom}[1]{[{#1}]_{K\!K}}
\newcommand{\KK}{{\mbox{$K\!K$}}}
\newcommand{\KOSdo}{\KKK(\Hdo)\otimes SB}
\newcommand{\KOS}{\KKK(\Hwh)\otimes SB}

\newcommand{\KOam}[1]{\KKK(H_{#1})\otimes B}
\newcommand{\KOdo}{\KKK(\Hdo)\otimes B}
\newcommand{\KOtwo}{\KKK(\Htwo)\otimes B}
\newcommand{\KOg}{\KKK(H)\otimes B}
\newcommand{\KOo}{\KKK(\Hha)\otimes B}
\newcommand{\KOquad}{\KKK(\Hquad)\otimes B}
\newcommand{\KOsix}{\KKK(H_{4}) \otimes B}
\newcommand{\KO}{\KKK(\Hwh)\otimes B}
\newcommand{\KVT}{\cite{ggk:hcmtst}}
\newcommand{\Kspone}{K_{*+1}}
\newcommand{\Ks}{K_*}
\newcommand{\ho}{homomorphism}
\newcommand{\hos}{homomorphisms}
\newcommand{\st}{such that   }
\newcommand{\LLL}{{\mathcal L}}
\newcommand{\MDrho}{D_\MZtau}
\newcommand{\MMSdo}{M(\KOSdo)}
\newcommand{\MMS}{M(\KOS)}
\newcommand{\MMam}[1]{M(\KOam{#1})}
\newcommand{\MMdo}{{M(\KOdo)}}
\newcommand{\MMg}{{M(\KOg)}}
\newcommand{\MMo}{M(\KOo)}
\newcommand{\MMquad}{{M(\KOquad)}}
\newcommand{\MM}{{M(\KO)}}

\newcommand{\MZtauco}{\CALrho(A)\COM}
\newcommand{\MZtaucos}{\CALrho_2(A)\COM}
\newcommand{\MZtau}{{{\Theta}}}

\newcommand{\Mtaua}{{{\gamma}}}
\newcommand{\Mtau}{{{\theta}}}
\newcommand{\Mus}{\Mu_s}
\newcommand{\Mut}{{\Mu_t}}
\newcommand{\Mu}{{{u}}}
\newcommand{\NN}{{\mathbb{N}}}
\newcommand{\Odo}{{1_{3}}}

\newcommand{\Pext}{\operatorname{Pext}}

\newcommand{\QQQ}{{\mathbb Q}}
\newcommand{\QQS}{Q(\KOS)}

\newcommand{\QQdo}{{Q(\KOdo)}}
\newcommand{\QQg}{{Q(\KOg)}}

\newcommand{\QQquad}{{Q(\KOquad)}}
\newcommand{\QQsix}{{Q(\KOsix)}}

\newcommand{\RR}{{\mathbb{R}}}
\newcommand{\UG}{{\mathcal U}}

\newcommand{\ZZp}[1]{\ZZ/{#1}}
\newcommand{\ZZ}{{\mathbb{Z}}}

\newcommand{\Zha}{{0_{1}}}

\newcommand{\abequiv}{\sim}
\newcommand{\ampliC}{\otimes C}

\newcommand{\atsign}{\char64\relax}

\newcommand{\bigproj}{\operatorname{Proj}}
\newcommand{\chif}{ \chi^{\phi}}

\newcommand{\cstar}{\mbox{$C^*$}}
\newcommand{\dist}{\operatorname{dist}}
\newcommand{\dotPhi}{\dot{\Phi}}
\newcommand{\dotPsi}{\dot{\Psi}}
\newcommand{\dotT}{\dot{T}}
\newcommand{\dotI}{\dot{I}}
\newcommand{\eep}{e^\perp}
\newcommand{\ee}{e}

\newcommand{\etaKs}{{\eta^{*}}}
\newcommand{\fseta}{{\mathcal G}}

\newcommand{\fset}{{\mathcal F}}
\newcommand{\hugePK}[1]{\underline{\mathbf K}\left({#1}\right)}
\newcommand{\huK}{\underline{\mathbf K}}
\newcommand{\id}{\operatorname{id}}
\newcommand{\infer}{\Longrightarrow}
\newcommand{\inv}{^{-1}}
\newcommand{\kone}[2]{\left[{#1}\right]_{#2}}
\newcommand{\lind}{d_\iota}
\newcommand{\matrM}{{\mathbf M}}
\newcommand{\mul}[1]{{#1}\!\cdot\!}
\newcommand{\smul}[1]{{#1}\cdot}
\newcommand{\myPrSu}{\left.\myPr\right/\mySu}
\newcommand{\myPr}{{\prod{B_i}}}
\newcommand{\mySu}{{\sum{B_i}}}
\newcommand{\mysmPrSu}{\left.\mysmPr\right/\mysmSu}
\newcommand{\mysmPr}{{\textstyle \prod{B_i}}}
\newcommand{\mysmSu}{{\textstyle \sum{B_i}}}

\newcommand{\norm}[1]{{\left\|{#1}\right\|}}

\newcommand{\pKstar}[2]{K_*({#2};\ZZp{#1})}
\newcommand{\phis}{\phi}
\newcommand{\phit}{\phi}
\newcommand{\projset}{{\mathcal P}}
\newcommand{\psis}{\psi}
\newcommand{\psit}{\psi}
 
\newcommand{\spectrum}{\operatorname{sp}}

\newcommand{\vvp}{v^\perp}
\newcommand{\vv}{v}
\renewcommand{\Im}{\operatorname{Im}}
\renewcommand{\epsilon}{\varepsilon}
\newcommand{\ep}{\epsilon}
\renewcommand{\phi}{\varphi}

\newcommand{\OOO}{{\mathcal O}}
\newcommand{\myphi}{\overline{\phi}}
\newcommand{\mypsi}{\overline{\psi}}

\newcommand{\twobytwo}[4]{\left[\begin{matrix}{#1}&{#2}\\{#3}&{#4}\end{matrix}\right]}
\newcommand{\finiteref}{(iii.1)}
\newcommand{\finitereftwo}{(iv.1)}
\newcommand{\infiniteref}{(iii.$\infty$)}
\newcommand{\diag}{\operatorname{Diag}}
\newcommand{\projsetKs}{{\mathcal P}_*}
\newcommand{\projsetnil}{{\mathcal P}_0}
\newcommand{\uniset}{{\mathcal V}}
\newcommand{\projsetKhu}{{\mathcal P}}
\newcommand{\pstar}{_{\sharp}}
\newcommand{\unPr}{1_{\Pi}}
\newcommand{\unPrSu}{1_{\Pi/\Sigma}}
\newcommand{\DD}[1]{{\mathbb I}^\sim_{#1}}
\newcommand{\DDn}[1]{{\mathbb I}_{#1}}
\newcommand{\SCHECK}{}
\newcommand{\SCOMMENT}[1]{}
\newcommand{\rera}{\mathbf{rr}}
\newcommand{\vali}{\mathbf{cco}}
\newcommand{\valii}{\mathbf{pfo}}
\newcommand{\valiii}{\mathbf{elo}_0}
\newcommand{\valiv}{\mathbf{elo}_1}
\newcommand{\valv}{\mathbf{ipo}}
\newcommand{\PRODrefuhETAinj}{Theorem \ref{collect}(i)}
\newcommand{\PRODrefuhHOMinj}{Theorem \ref{collect}(ii)}
\newcommand{\PRODrefuhCORO}{Theorem \ref{collect}(iv)}
\newcommand{\PRODrefuhPext}{Theorem \ref{collect}(iii)}
\newcommand{\Aff}{\operatorname{Aff}}
\begin{document}
\title{On the classification of nuclear
 \cstar-algebras}
\author{Marius Dadarlat\thanks{Partially supported by the National
 Science Foundation, \# 9622434-DMS} \and S\o ren
 Eilers\thanks{Partially supported by the Carlsberg Foundation}}
\date{Final version, August 1998}
\maketitle
\tableofcontents

\section{Introduction}
Two of the most influential works on \cstar-algebras from the
mid-seventies --- Brown,  Douglas and Fillmore's \cite{bdf:eck} and Elliott's \cite{gae:cilssfa} --- both contain 
\emph{uniqueness} and \emph{existence} results in the now standard sense which we
shall outline below.
These papers served as keystones for two separate 
theories ---\emph{$\KK$-theory} and \emph{the classification program}  --- which for many years
parted ways with only moderate interaction. But with this common
origin in mind, it is not surprising that recent
years have seen a fruitful interaction which has been one of the main
engines behind rapid progress in the classification program.

In the present paper we take this interaction even further. 
We prove general existence and uniqueness
results using $\KK$-theory and a concept of quasidiagonality for
representations. These results are employed to obtain new
classification results for certain classes of quasidiagonal
\cstar-algebras introduced by H.\ Lin.
An important novel feature of these classes is that they are
defined by a certain
local approximation property, rather than by
 an inductive limit construction.

Our existence and uniqueness results are in the spirit of
the classical $\Ext$-theory from \cite{bdf:eck}. The main
complication overcome in the paper is to control the stabilization
which is necessary when one works with finite \cstar-algebras. In the
infinite case, where programs of this type have already been 
successfully carried out, stabilization is unnecessary. Yet, our methods are
sufficiently versatile to allow us to reprove, from a handful of basic
results, the classification of purely infinite nuclear \cstar-algebras
of Kirchberg and Phillips.

Indeed, it is our hope that
this can be the starting point of a unified approach to classification
of nuclear \cstar-algebras.

Apart from $\KK$-theory, the main technical tools are approximate
morphisms and partial $\KK$-elements, defined using $K$-theory with
coefficients. 

\subsection{Existence and uniqueness}

Existence and uniqueness theorems can be found in most, if not all,
classification papers to 
combine with refinements of Elliott's intertwining argument
to yield classification theorems. 

In this framework, a \textbf{uniqueness result} allows one to conclude that if
the $K$-theoretical elements induced by two $*$-homomorphisms
$\phi,\psi:A\longrightarrow B$ coincide, with $A,B$ \cstar-algebras
satisfying extra conditions, then $\phi$ and $\psi$ are
equivalent in a sense which generalizes unitary equivalence.
 In \cite{gae:cilssfa}, one concludes from
$\phi_*=\psi_*$ on $K_0(A)$ that $\phi$ and $\psi$ are approximately
unitarily equivalent when $A$ and $B$ are both $AF$ (cf.\ \cite{ob:af}). In \cite{bdf:eck} one
shows that if
 $\KKhom{\phi}=\KKhom{\psi}$,  then
$\phi$ and $\psi$ are unitarily equivalent whenever
$A=C(X)$ and $B=\LLL(H)/\KKK(H)$.

The role of the equally important \textbf{existence results} is to
provide means of realizing given $K$-theoretical data by a sufficiently
multiplicative \CPAM. In \cite{gae:cilssfa}, one realizes any positive element of
$\Hom(K_0(A),K_0(B))$ by a $*$-homomorphisms when $A$ and $B$ are
$AF$.  In \cite{bdf:eck}, any element of $\KK(C(X),\LLL(H)/\KKK(H))$ is
realized by a $*$-monomorphism.

One of the main obstacles in achieving general existence and
uniqueness results has been the fact that,  as soon as one ventures
beyond these classical examples, one only can expect to achieve
\emph{stable} versions of such theorems in a sense that involves
adding or subtracting maps of the form $\mu:A\longrightarrow \matrM_n(B)$.
The challenge has been to control the complexity of the stabilization,
in order to be able to  incorporate them into classification results.

This has been quite successfully carried out, using maps $\mu$ with
finite-dimensional images, in important classes of quasidiagonal
\cstar-algebras. But 
in most existence and uniqueness results so far  the source $A$ has been required to
be a member of a small class of \cstar-algebras forming \emph{building
blocks} for the class of \cstar-algebras one has tried to classify,
leading to restrictions on the ensuing classification results. 

We offer, in the present paper, existence and uniqueness results valid
for sources way beyond even the full class of nuclear quasidiagonal
\cstar-algebras. More precisely, we consider a
unital and nuclear source $A$ and a unital target $B$, and require
that this pair allows an \emph{absorbing}  and
\emph{quasidiagonal} unital representation $\gamma:A\longrightarrow \MMg$ 
as defined in Section \ref{uniqprelims} below. Such a map exists
automatically when $A$ is quasidiagonal, or (by a theorem of Lin \cite{hl:saueh}) when $A$ can be embedded
into $B$ via a third simple \cstar-algebra.

In this case, our uniqueness result Theorem \ref{vVaumPROTO} states 
that if $\phi,\psi: A \longrightarrow B$
are two $*$-\hos\ inducing the same \KK-class, we may conclude that $\phi$ 
is stably approximately unitarily equivalent to $\psi$ in a sense
involving adding ``finite pieces'' of $\gamma$.
For \cstar-algebras $A$ which satisfy the universal coefficient
theorem of \cite{jrcs:ktuctkgk}, our result is predated by one of a similar
nature, valid when $A$ or $B$ is simple, which appears in \cite{hl:saueh}. We emphasize that our
uniqueness result does
not depend on the universal coefficient theorem, and we believe that
in addition to proving a more general result, our arguments are
somewhat more
conceptual.

In fact, we shall require uniqueness results which hold also for
sufficiently multiplicative \CPAM
s, and this entails the problem of associating 
to such maps a kind of  \emph{partial} \KK-elements to substitute for the
globally defined group homomorphisms one gets from
$*$-homomorphisms. This is 
done using the universal multicoefficient theorem of
\cite{mdtal:umctkg} in a fashion explained by  the first named author
in his talk
at the \emph{Workshop on the classification of amenable \cstar-algebras} at
the Fields Institute in December of 1994. As soon as we have made
sense out of this concept, a uniqueness result for \CPAM s, Theorem
\ref{vVII}, can be derived from
those for $*$-homomorphisms using a procedure originating with \cite{hlncp:ammco} in  the  
torsion free case. This method also requires keeping a close eye on
the $K$-theory for products of \cstar-algebras.

In our existence result Theorem \ref{EvI} we manage to realize -- 
partially, in a sense corresponding to the one used in the uniqueness
result -- a given element 
from $\KK(A,B)$ as a difference of 
\CPAM s from $A$ to $\matrM_N(B)$. Again, all we require is that $A$
is nuclear, and that  an absorbing quasidiagonal representation $\gamma:A\longrightarrow\MMg$
exists as in the
uniqueness case. Furthermore, one of the maps in the difference can be
chosen as a ``finite piece'' of $\gamma$.

In the building block approach,
establishing existence has typically been somewhat easier than achieving
uniqueness results. This may still be the case in our setting, but at
the current stage it is in fact the
existence which is causing problems. Indeed, the existence result
offered in our paper has shortcomings in the finite case which are
responsible for a number of unwanted, and hopefully redundant,
restrictions in the resulting classification results.
For instance, we do not have sufficient technology to prove in full
generality that a
positive \KK-element can be 
realized by a single map, as one would expect to be the case.

As a main application, we apply our existence and uniqueness results
to the class of TAF \cstar-algebras introduced and studied by H.\
Lin. Using a factorization property of these \cstar-algebras, we
are able to prove in Theorem \ref{reduction} that, up to an isomorphism, there is only one 
unital, separable, nuclear and simple TAF \cstar-algebra satisfying the UCT with
$K_0(A)=\QQQ$ and $K_1(A)=G$, where $G$ is a countable fixed arbitrary group, thus proving that every such \cstar-algebra
falls in the well studied class of $AD$ algebras of real rank zero.

We believe, however, that our existence and uniqueness results will be
applicable to a wide range of classification problems, extending well
beyond the TAF case. In fact, it might be that at least the
uniqueness result will be sufficiently  versatile to serve as a
unifying element for many future efforts to classify nuclear
\cstar-algebras. To substantiate this claim, we apply our methods to
the case of purely infinite \cstar-algebras, reproving in Theorem
\ref{pisunclas}, rather easily,
the classification theorem of purely infinite simple
unital nuclear \cstar-algebras (see \cite{ek:cpicukt} and \cite{ncp:ctnpisc})
from a handful of fundamental results about such algebras. 

We reported on the present paper at the 1998 GPOTS. At the same
conference, H.\ Lin reported results which -- although the methods
differ -- overlap with our classification results in the  TAF
\cstar-algebra case. More details are given in the notes of the present paper.
\subsection{Methods}
To arrive at uniqueness and existence results from the existence of an
absorbing and quasidiagonal representation, we depend on the full force of \KK-theory,
including several of the different
realizations of the Kasparov groups
and their interrelations. 

To achieve such representations, one may use  the Kasparov Weyl-von
Neumann-Voiculescu theorem along with the concept of
quasidiagonality. One can also employ
a related result by H.\ Lin
\cite{hl:saueh} stating that a certain extension associated to a
unital inclusion $\iota:A\longrightarrow B$ is absorbing when $A$ is
nuclear and  $B$ is
simple. 

To refine uniqueness results we depend on a number of basic properties
about the $K$-theory of products of \cstar-algebras. Furthermore, to
apply our results to classify  TAF \cstar-algebras,
we use several results about structural properties of such
\cstar-algebras by Lin. We prove classification for purely infinite
algebras by appealing to the embedding theorem for exact
\cstar-algebras of  \cite{ekncp:eeccfco}, as well as structural
results by Kirchberg, Phillips and R\o rdam.
\subsection{Organization}
The paper is organized as follows. In Section \ref{gepr} we lay out
notation and define several relevant classes of \cstar-algebras. Then
in Section \ref{abunique}, we establish (without using the universal coefficient theorem) the uniqueness result which is
at the core of the paper. This is subsequently
refined (using, among other things, the universal coefficient theorem)
in Section \ref{doctor}. The basic existence results are collected in
Section \ref{exist}. The applications to classification are presented
in Sections \ref{tafclas} and \ref{piclas}, concerning the TAF and the
purely infinite \cstar-algebras, respectively. Finally, Appendix
\ref{KKap} contains results about the $K$-theory of products of
\cstar-algebras and an explanation of how to associate partial
$\KK$-elements to sufficiently multiplicative completely positive
maps.

\section{General preliminaries}\label{gepr}
\subsection{Notation}

\subsubsection*{Some classes of \cstar-algebras}

We single out several classes of \cstar-algebras for easy reference:

\begin{defin}{UCTalg}
We say that a separable \cstar-algebra \emph{satisfies the UCT} if the diagram
\[
\xymatrix{
0\ar[r]&{\Ext(\Ks(A),\Kspone(B))}\ar[r]&
{\KK(A,B)}\ar[r]&{\Hom(\Ks(A),\Ks(B))}\ar[r]&0
}
\]
is a short exact sequence for every $\sigma$-unital algebra $B$.
\end{defin}

A large class of algebras satisfying UCT was exhibited in
\cite{jrcs:ktuctkgk}. It is not known whether there exist separable nuclear
\cstar-algebras not satisfying the UCT.
If the separable \cstar-algebra $A$  satisfies the UCT, then
for any $\sigma$-unital \cstar-algebra $B$ the sequence:
\begin{equation} \label{UMCT}
\xymatrix{
0\ar[r]&{\Pext(\Ks(A),\Kspone(B))}\ar[r]&
{\KK(A,B)}\ar[r]&{\Hom_\Lambda(\hugePK{A},\hugePK{B})}\ar[r]&0
}
\end{equation}
is also exact by
 \cite{mdtal:ccomk}. Here  $\hugePK{-}$ denotes the sum of all $K$-theory groups with $\mathbb {Z}/n$
coefficients, $n \geq 1$, and $\Lambda$ denotes the natural set of coefficient
transformations and the Bockstein operations (see \cite{cs:tmcIV} and \cite{mdtal:ccomk}).

\begin{defin}{simpleinclu}
We say that a unital $*$-homomorphism $\iota:A\longrightarrow B$ is a
\emph{unital simple embedding} if it can be factored 
\[
\xymatrix{
{A}\ar@{^(->}[r]^{\iota'}&{C}\ar@{^(->}[r]^{\iota''}&{B}
}
\]
where $\iota',\iota''$ are injective and $C$ is simple.
\end{defin}
This implies that $C, \iota'$ and $\iota''$ are unital. Note that the composition (on either side) of a
simple embedding with a unital injective  $*$-homomorphism is a simple  embedding.


%
%
%
%
If $B$ is a C*-algebra we denote by $\bigproj(B)$ the set of all selfadjoint
projections in $B$.
The K-theory class of a projection $p$ is denoted by $[p] \in K_0(B)$.
If $B$ is unital we let $\UG_{n}(B)$ denote the unitary group of $\matrM_n(B)$.

\begin{defin}{targetREST}
A \cstar-algebra $B$ is called an \emph{admissible target algebra}
if it is unital, has real rank zero (\cite{lgbgkp:crrz}) and satisfies
\begin{rolist}
\item whenever $p,q\in \bigproj(B\otimes\KKK)$, then $[p]=[q]\Longrightarrow p\oplus 1_B\sim q\oplus
1_B$
\item the canonical map $\UG_{1}(B)\longrightarrow K_1(B)$ is surjective
\end{rolist}
and if either
\begin{list}{}{\setlength{\leftmargin}{0.99cm}}
\item[\infiniteref] the canonical map $\bigproj(B)\longrightarrow
K_0(B)$ is surjective
\end{list}
or both of
\begin{list}{}{\setlength{\leftmargin}{0.99cm}\setlength{\parsep}{0cm}}
\item[\finiteref] For any $x \in K_0(B)$ such that $nx \geq 0$ for some $n \geq 0$,
one has $x+[1_B] \geq 0$.
\item[\finitereftwo]  For any $x\in K_0(B)$  and any $n\not=0$, there is
$y\in K_0(B)$ such that $-[1_B]\leq y\leq [1_B]$ and $x-y\in nK_0(B)$.
\end{list}
holds.
\end{defin}

When needed, we distinguish between admissible target algebras satisfying
\finiteref--\finitereftwo\ or  \infiniteref\ by calling them admissible of \emph{
finite type} or \emph{infinite type}, respectively. Examples will be
given in Propositions \ref{isTAFtarget} and \ref{isINFTYtarget}.

The point of this definition is that whenever a sequence of admissible
targets are given, then both of the natural maps
\begin{gather}
\Hom_\Lambda(\hugePK{A},\hugePK{\mysmPr})
\longrightarrow \prod{\Hom_\Lambda(\hugePK{A},\hugePK{B_i})}\label{injone}\\
\KK(A,\mysmPrSu)\longrightarrow
\Hom_\Lambda(\hugePK{A},\hugePK{\mysmPrSu})\label{injtwo}
\end{gather}
will be injective (the latter in fact an isomorphism)
 for an $A$ satisfying the UCT.
We defer the proof of this to Appendix \ref{KKinj} below.
We are also going to need that 
whenever $B_i$ is a sequence of admissible target
algebras, then $\myPrSu$ is admissible.

\begin{remar}{less}
In fact, we can get injectivity for the maps discussed above asserting
considerably weaker versions of (iii)-(iv).
 This will be clear from
Appendix \ref{KKinj}. On the other hand, $B_i=C([0,1])$ is a
counterexample to injectivity in \eqref{injone} and $B_i={\mathbb C}$ a counterexample to \eqref{injtwo}.
\end{remar}

\subsection{\KK-theory}\label{KKpictures}

We depend on Kasparov's \KK-theory from \cite{ggk:okec} for most of
this paper, as well on several different characterizations or
realizations of it. As a
standard picture of a \KK-cycle, we shall adopt the one from
\cite[2.1]{nh:ck}. There $\KK(A,B)$ 
is defined in terms of triples $(\phi_+,\phi_-,x)$ where
$\phi_\pm:A\longrightarrow \MMg$ are $*$-homomorphisms and $x\in\MMg$
satisfies
\begin{gather}\label{KKone}
x\phi_+(a)-\phi_-(a)x\in\KOg\\ \label{KKtwo}
\phi_+(a)(x^*x-1),\phi_-(a)(xx^*-1)\in\KOg
\end{gather}
for each $a\in A$. Higson works with separable \cstar-algebras only,
but his picture extends readily to the case of a $\sigma$-unital
$B$. Equally important to us is the related Cuntz picture of \KK-theory,
studied in \cite{nh:tt},
 which we
shall indicate by $\KK_h$ as in Chapter 4 of \cite{kkjkt:ek}. Here
the cycles are \emph{Cuntz pairs} $(\phi_+,\phi_-)$ satisfying
\[
\phi_+(a)-\psi_-(a)\in\KOg.
\]

\section{Uniqueness up to absorption}\label{abunique}
In the first section of the paper we establish, via $\Ext$-
and $\KK$-theory, a uniqueness result proving that two $\KK$-equivalent
\CPAM s are approximately unitarily equivalent after one adds an
absorbing representation. Most of the work goes into controlling the nature
of the unitary implementing this equivalence, proving that it interacts
sufficiently well with the other components to allow truncation in a
sense to be made precise in Section \ref{doctorPROTO}.
\subsection{Preliminaries}\label{uniqprelims}
\subsubsection*{Notation and conventions}
In all of Section \ref{abunique}, we only work with infinite, separable Hilbert spaces, so all Hilbert
spaces in this paper are isomorphic. However, we introduce the following
notation to aid the reader in distinguishing between different
instances of them. We start with a separable
Hilbert space $\Hha$ and define
\[
H_m=\overbrace{\Hha\oplus\dots\oplus \Hha}^m.
\]
for any $m\in\NN$. There are now
canonical identifications between, say, $\matrM_2(\KOo)$ and $\KOtwo$,
and we shall employ them tacitly in the following. However,  we choose
\emph{not} to apply  the (non-canonical) isomorphisms between, e.g.,
$\KOo$ and $\KOtwo$, as we feel
this helps to clarify our constructions. We shall abandon this
practice towards the end of this section.

We work with the multiplier
algebras $M(\KKK(H_m)\otimes B)$, as well as the corona algebras
\[
Q(\KKK(H_m)\otimes B)=M(\KKK(H_m)\otimes B)/\KKK(H_m)\otimes B.
\]
The quotient map from $M(\KKK(H_m)\otimes B)$ to $Q(\KKK(H_m)\otimes
B)$ is denoted by $\pi_m$, and whenever there is need for distinction,
we write $1_m$ and $0_m$ for the
identity and zero elements of these algebras. 

\begin{defin}{tauconv}
An \emph{admissible scalar representation}
$
\Mtau:A\longrightarrow\MMo
$
is a $*$-homomorphism which factors as
\[
\xymatrix{
{A}\ar[r]^-{\Mtau'}&{\LLL(H)}\ar[r]^-{-\otimes\id}&{{\mathcal
L}(H)\otimes M(B)}\ar@{^(->}[r]&{\MMg}
}
\]
where $\Mtau'$ is unital, faithful and of infinite multiplicity, i.e.\ of the
form $\mul{\infty}\Mtaua$ for some representation $\Mtaua$. 
\end{defin}

When $\Mtau$ is a admissible scalar representation, we are also going to
consider representations of the form
\[
0_{m-1}\oplus\Mtau:A\longrightarrow \MMam{m}
\]
for $m \geq 2$. When the size of
the added zero is clear from the context, or irrelevant, we denote this
representation by $\MZtau$. Note that by convention, $\MZtau$ is \emph{never}
unital.
\subsubsection*{Absorbing and quasidiagonal representations}
\begin{defin}{abequivalent}
Fix a unital \cstar-algebra $B$. When $\Mtaua:A\longrightarrow \LLL_B(E)$ and
$\Mtaua':A\longrightarrow 
\LLL_B(E')$ are two unital representations, with $E$ and $E'$ Hilbert
\cstar-modules over $B$, we say that $\Mtaua$ and $\Mtaua'$ are
\emph{equivalent}, and write $\Mtaua\abequiv\Mtaua'$, if there exists
a sequence $U_m\in\LLL_B(E,E')$, consisting of unitaries, 
such that
\begin{rolist}
\item $\norm{\Mtaua(a)-U_m\Mtaua'(a)U_m^*}\longrightarrow0$,
$m\longrightarrow\infty$
\item $\Mtaua(a)-U_m\Mtaua'(a)U_m^*\in\KKK_B(E)$
\end{rolist}
for any $a\in A$.
\end{defin}
\begin{defin}{absorption}
A unital representation $\Mtaua:A\longrightarrow \LLL_B(E)$ is
\emph{absorbing} if for any other unital representation
$\Mtaua':A\longrightarrow \LLL_B(E')$,
$\Mtaua\oplus\Mtaua'\abequiv\Mtaua$. 
\end{defin}

Clearly any two absorbing representations are equivalent.
Kasparov proved in \KVT,
generalizing Voiculescu's result in \cite{dv:nwvt}, that any admissible
scalar representation is absorbing if $A$ is separable and nuclear.
Another class of absorbing representations was exhibited by Lin in \cite{hl:saueh}, based on unital
inclusions of $A$ into $B$, where either $A$ or $B$ is simple. The
observation that Lin's proof carries over to the case of unital simple  embeddings
(cf. \ref{simpleinclu}) is so crucial to our approach that we shall
state is as a separate lemma:

\begin{lemma}{siabsorbs}
When $\iota:A\longrightarrow B$ is a unital simple  embedding
and $A$ is nuclear,
the map $\lind:A\longrightarrow \MMg$ given by
$\lind(a)=1\otimes\iota(a)$
is absorbing.
\end{lemma}
\begin{demo}
By assumption, a simple \cstar-algebra $C$ and unital inclusions $\iota':A\longrightarrow C$ and
$\iota'':C\longrightarrow B$ can be found with $\iota=\iota''\iota'$. By \cite[1.6]{hl:saueh},
$d_{\iota'}$ is absorbing. Since  $\id\otimes\iota''$ preserves
approximate units, it induces a map
$\widehat{\iota''}:M(\KKK\otimes C)\longrightarrow M(\KKK\otimes B)$, and by
\cite[1.11]{hl:saueh} (valid since $A$ is nuclear), we conclude that
$\widehat{\iota''}d_{\iota'}=d_{\iota''\iota'}=\lind$ is also absorbing.
\end{demo}

\begin{defin}{quasidiagonalrep}
A representation $\Mtaua:A\longrightarrow \MMg$ is
\emph{quasidiagonal} if there exists an approximate unit of projections
$(e_n)$ for $\KOg$ with the property that
\[
[e_n,\Mtaua(a)]\longrightarrow 0
\]
for all $a\in A$.
\end{defin}

If $\Mtaua$ is quasidiagonal and  $\Mtaua \sim  \Mtaua'$, then
 $\Mtaua'$ is quasidiagonal.
It is clear that if $\Mtaua$ is quasidiagonal, one may
assume that $(e_n)$ has the property that $e_n\in \matrM_{r_n}(B)$ 
for some sequence of integers $(r_n)$. 

When this is the case, we may define a sequence of completely positive
morphisms $\Mtaua_n:A\longrightarrow \matrM_{r_n}(B)$  by
$\Mtaua_n(a)=e_n\Mtaua(a)e_n$. We call $(\Mtaua_n)_{n\in\NN}$ a
\emph{quasidiagonalization} of $\Mtaua$ by $(e_n)_{n\in\NN}$.
Note that the sequence  $(\Mtaua_n)_{n\in\NN}$ is asymptotically
multiplicative.

The concept of relative quasidiagonality that we use in this paper
was studied in \cite{ns:rqaKKt}
 and \cite{cs:fskg}.
\begin{remar}{absqdnotes}
If there exists an absorbing and quasidiagonal representation
$\gamma:A\longrightarrow\MMg$, then \emph{any} absorbing
representation will be quasidiagonal. 

Note that $\lind$ of Lemma \ref{siabsorbs} is always
quasidiagonal, as it commutes with projections $e_n=\mul{n}1_B$. Thus all absorbing representations
$\gamma:A\longrightarrow \MMg$ are quasidiagonal when there is a
unital simple embedding of $A$ into $B$ and $A$ is nuclear. The same holds true whenever
$A$ is a nuclear quasidiagonal \cstar-algebra.
\end{remar}

\subsection{Consequences of $\Ext$-theory}
\subsubsection*{Obtaining elements of $\Ext$}

We work with \hoos\  $\phi:A\to \KO$. 
If $\phi :A \to B$ is a \hoo, then we regard $\phi$ as a map into
$\KO$ by embedding $B$ as a $(1,1)$-corner of $\KO$.
The pullback of the essential, semisplit extension
\begin{equation} \label{C}
0 \to \KOS  \to \KKK(\Hwh) \otimes CB \to \KO \to 0
\end{equation}
by the \hoo\ $\phi$ is the mapping cone extension
\begin{equation} \label{mcone}
0 \to \KOS  \to C_\phi \to A \to 0
\end{equation}
This is an essential, semisplit extension. 
With
\[
\chif(a)(t)=t\phit(a), \quad a \in A,\,\, t \in [0,1],
\]
we get a completely positive\SCHECK{} contractive map $\chif:A\rightarrow
C_b([0,1),\KO)$. Since 
\[
C_b([0,1),\KO)\subseteq C_{b,\operatorname{strict}}((0,1),\MM)=\MMS,
\]
we may and shall consider $\chif$ as a map into $\MMS$. 
With $\chif$ as above, one checks that
 $\pi_1\circ\chif:A\longrightarrow \QQS$ is
the Busby invariant of (\ref{mcone}).
Since $\chif$ is completely positive and contractive, we have that 
(\ref{mcone}) is semisplit as claimed, and hence it defines an element
$[C_\phi] \in \Ext(A,SB)\inv$.
\begin{propo}{extequal}
Let $A$, $B$ be \cstar-algebras with $A$ separable and $B$ $\sigma$-unital.
\newline Let $\phi,\psi: A \to \KO$ be two \hoos.
If $[\phi]_{KK}=[\psi]_{KK}$ in $KK(A,B)$, then $[C_\phi]=[C_\psi]$
in $\Ext(A,SB)\inv$.
\end{propo}
\begin{demo}
Since the isomorphism $\gamma:\Ext(X,Y)\inv \to KK^1(X,Y)$ of Kasparov
is natural, we have a commutative diagram
\[\xymatrix{
\Ext(B,SB)\inv \ar[r]^{\phi^*} \ar[d]_{\gamma}&\Ext(A,SB)\inv \ar[d]^{\gamma}\\
KK^1(B,SB)\ar[r]^{\phi^*}&KK^1(A,SB) }\]
and a similar diagram for $\psi$.
If $x$ is the class of (\ref{C}) in $\Ext(B,SB)\inv$, then $[C_\phi]=\phi^*(x)$
and  $[C_\psi]=\psi^*(x)$.
Since $\gamma$ is an isomorphism, in order to show that
$[C_\phi]=[C_\psi]$ it suffices to prove
$\gamma\phi^*(x)=\gamma\psi^*(x)$. Using the commutative diagram above,
this is equivalent to showing that $\phi^*\gamma(x)=\psi^*\gamma(x)$.
By \cite[18.7.2]{bb:koa}, for any $y \in KK^1(B,SB)$, $\phi^*(y)$ 
equals the Kasparov product $[\phi]_{KK} \otimes y$.
Therefore
\[ \phi^*\gamma(x)=[\phi]_{KK} \otimes \gamma(x)=
[\psi]_{KK} \otimes \gamma(x)= \psi^*\gamma(x). \]
\end{demo}

\subsubsection*{Weak uniqueness}

The starting point of our investigation of uniqueness of maps between two
\cstar-algebras $A$ 
and $B$ will be two \hoos, $\phi$ and $\psi$. We will require $A$ to
be separable, but it is crucial for applications in Section
\ref{doctor} that $B$ can be any $\sigma$-unital \cstar-algebra.

\begin{propo}{vI}
Let $A$ be a unital, separable, nuclear C*-algebra and let $B$ be a  $\sigma$-unital \cstar-algebra.
Assume that $\phi$ and $\psi$ are two \hoos\ from $A$ to $\KO$ which satisfy
$[\phi]_{KK}=[\psi]_{KK}$ in $KK(A,B)$.
Then for any admissible scalar
representation $\Mtau:A\rightarrow\MMo$, there exists
a strictly continuous map
\[
\Mu:(0,1)\longrightarrow\UG(\MMdo)
\]
with the property that
\[
\Mut\left[\begin{matrix}t\phi_t(a)&&\\&\Zha&\\&&\Mtau(a)\end{matrix}\right]
\Mut^*-\left[\begin{matrix}t\psi_t(a)&&\\&\Zha&\\&&\Mtau(a)\end{matrix}\right]\in
C_0((0,1),\KOdo)
\]
for all $a\in A$.
\end{propo}
\begin{demo}
By Proposition \ref{extequal} 
we conclude that $\pi_1 \circ\chi^\phi$ and $\pi_1\circ\chi^\psi$ define the same
element of $\Ext(A,\KOS)$. Since $0_1\oplus\Mtau$ defines an
absorbing extension by Kasparov's
Voiculescu theorem \KVT, 
we have that
$\pi_3\circ(\chi^\phi\oplus\Zha\oplus\Mtau)$ is equivalent to 
$\pi_3\circ(\chi^\psi\oplus\Zha\oplus\Mtau)$. This means that
\[
\Ad_{\pi_3(u)}\circ\pi_3\circ\left(\chi^\phi\oplus\Zha\oplus\Mtau\right)=
\pi_3\circ\left(\chi^\psi\oplus\Zha\oplus\Mtau\right)
\]
for some $u\in\UG(\MMSdo)$. In $\MMSdo$, this amounts
exactly to saying
that the difference considered in the Proposition is an element of
$\KOSdo$. 
\end{demo}

%

\subsubsection*{A trivial \KK-cycle}
Working with $\Mut$ and $\Mtau$ as given in Proposition \ref{vI}, we are going to consider
\[
\chi_0(a)=\MZtau(a)=\left[\begin{matrix}0_2&\\&\Mtau(a)\end{matrix}\right]
\qquad
\chi_t(a)=\Mut\left[\begin{matrix}0_2&\\&\Mtau(a)\end{matrix}\right]\Mut^*,
\]
where  $t\in (0,1)$. This way
we get a family of $*$-homomorphisms
\[
\chi_t:A\rightarrow \MMdo, \quad t \in [0,1)
\]
whose properties are collected in the following Lemma.

\begin{lemma}{phitprops}
Let  $\Mut$ and $\Mtau$ be as in Proposition \ref{vI}, and fix $t\in(0,1)$. For any $a\in A$, we have
\begin{rolist}
\item $\forall s\in [0,t]: \chi_0(a)-\chi_s(a)\in\KOdo$
\item $\chi_0(a)-\chi_s(a)$, considered as a
function from $[0,t]$ to $\KOdo$, is norm continuous.
\item $\chi_s(a)$, considered as a
function from $[0,t]$ to $\MMdo$,  is strictly continuous.
\end{rolist}
\end{lemma}
\begin{demo}
We have, with $\phi$ and $\psi$ as in Proposition \ref{vI},
\begin{eqnarray*}
\chi_0(a)-\chi_s(a)&=&
\left[\begin{smallmatrix}s\psis(a)&&\\&\Zha
&\\&&\Mtau(a)\end{smallmatrix}\right]
-
\left[\begin{smallmatrix}s\psis(a)&\\&0_2\end{smallmatrix}\right]-\\
&&
\Mus\left[\begin{smallmatrix}s\phis(a)&&\\&\Zha
&\\&&\Mtau(a)\end{smallmatrix}\right]\Mus^*
+
\Mus\left[\begin{smallmatrix}s\phis(a)&\\&0_2
\end{smallmatrix}\right]\Mus^*\\
&=&
\Mus\left[\begin{smallmatrix}s\phis(a)&\\&0_2
\end{smallmatrix}\right]\Mus^*-
\left[\begin{smallmatrix}s\psis(a)&\\&0_2
\end{smallmatrix}\right]+R_s(a),
\end{eqnarray*}
where $R_s(a)\in C_0((0,1),\KOdo)$ by Proposition \ref{vI}. Since the first two terms in the
last expression lie in $C_0((0,t],\KOdo)$, (i) and (ii) follow. That
(iii) holds follows from 
the fact that $s\mapsto \Mus$ is strictly continuous.
\end{demo}

\begin{propo}{trivkk}
Let  $\Mut$ and $\Mtau$ be as in Proposition  \ref{vI}, with $\MZtau=0_2 \oplus \Mtau$.
For any  fixed
$t\in(0,1)$,  $[\MZtau,\MZtau,\Mut]$ defines a trivial element
of $\KK(A,B)$.
\end{propo}
\begin{demo}
Fix $t\in (0,1)$. We first note, comparing Lemma \ref{phitprops}(i) to
\cite[4.1.1]{kkjkt:ek},
that  $(\chi_0,\chi_t)$ defines a 
cycle in
$\KK_h(A,B)$ (cf. Section \ref{KKpictures}). In fact, we get from (ii) and (iii) of Lemma
\ref{phitprops} that $(\chi_0,\chi_s)_{0\leq s\leq t}$
forms a homotopy, in the sense of \cite[4.1.2]{kkjkt:ek}, from
$(\chi_0,\chi_0)$ to $(\chi_0,\chi_t)$.
Since $[\chi_0,\chi_0]=0$,
 by \cite[4.1.4]{kkjkt:ek} 
 we conclude that
$[\chi_0,\chi_t]=0$ in $\KK_h(A,B)$. Applying the 
isomorphism
$ \mu_h:\KK_h(A,B)\rightarrow \KK(A,B)$
considered in \cite[4.1.8]{kkjkt:ek}, we get that $[\chi_t,\chi_0,1]$
is trivial in $\KK(A,B)$. Finally we note that, as explained for instance in \cite[2.3]{nh:cfek},
\[
[\chi_t,\chi_0,1]=
[\Ad_\Mut\circ\chi_0,\chi_0,1]=
[\chi_0,\chi_0,\Mut]=[\MZtau,\MZtau,u_t].
\]
\end{demo}
\subsection{Uniqueness up to absorption}
The preceding section left us with a trivial \KK-element
$[\MZtau,\MZtau,u]$ where $u$ was a unitary. Using results by Skandalis we now conclude that in
a very specific and rather subtle sense, the $K_1$-class of $u$ is
also trivial. The triviality translates to the fact that $u$ induces
an approximately inner automorphism on an auxiliary \cstar-algebra
defined by Lin.
\subsubsection*{A trivial $K_1$-element}
In this section we work with 
$\MZtau=0_2\oplus\Mtau:A\rightarrow \LLL_B(\Hdo\otimes B)=\MMdo$,
where $\Mtau$ is an
admissible scalar representation.
Recall that $A$ is unital and $\theta(1)=1$. We use the dot as a shorthand to indicate  composition by
the quotient maps $\pi:\MMg\longrightarrow \QQg$. Thus,
$\CALrho=\pi_3\circ\MZtau$ maps from $A$ to $\QQdo$. Furthermore, if 
$X\subseteq \QQg$, we denote by $X\COM$ the commutator of $X$
in $\QQg$.

We define a \cstar-algebra $\Drho$ by 
\[
\Drho=\{b\in \MMdo\mid [b,\MZtau(A)]\subset\KOdo\}.
\]
One checks directly that 
\begin{eqnarray}\label{Drhotwo}
\xymatrix{
{0}\ar[r]&
{\KOdo}\ar[r]^-{j}&
{\Drho}\ar[r]^-{\pi_3}&
{\CALrho(A)\COM}\ar[r]&0}
\end{eqnarray}
is a short exact sequence of \cstar-algebras.

\begin{propo}{inKoneimage}
With $A$ and $\MZtau$ as above, assume that $[\MZtau,\MZtau,u]=0$ in $\KK(A,B)$, where $u\in\MMdo$ is a
unitary. Then $u\in\Drho$, and
\[
\kone{u}{K_1(\Drho)}\in j_*(K_1(\KOdo)),
\]
where $j$ is the inclusion in \eqref{Drhotwo}. 
\end{propo}
\begin{demo}
Being part of a \KK-cycle, $u$ must commute with $\MZtau$ modulo the
compacts, and hence $u\in\Drho$.
We note that we have set up our \KK-cycle to be covered by the
description of 
$\KK(A,B)$ given in Proposition 2.6 of
\cite{gs:nnk}, where now $A$ is considered as a trivially
graded \cstar-algebra. Indeed, since our $\MZtau$ can
substitute as Skandalis' $\pi \otimes 1$, we may conclude from
$
[\MZtau,\MZtau,u]=[\MZtau,\MZtau,1_3]
$
that, by  Proposition 2.6 of \cite{gs:nnk},
\[
[\MZtau\oplus\MZtau,\MZtau\oplus\MZtau,u\oplus 1_3]\sim_{oh}
[\MZtau\oplus\MZtau,\MZtau\oplus\MZtau,1_3\oplus 1_3],
\]
where $\sim_{oh}$ denotes operator homotopy. 
This means that there is a norm continuous path of operators $\omega_s\in\MMquad$, $s \in [0,1]$ with
\[
\xymatrix{{u\oplus 1_3}\ar@{~}[rr]^-{\omega_s}&&{1_3\oplus 1_3}}.
\]
which satisfies, with $\MZtau_2=\MZtau\oplus \MZtau$
\begin{gather}\label{KOcommute}
[\MZtau_2(a),\omega_s]\in\KOquad\\
\MZtau_2(a)(\omega_s\omega_s^*-1),
\MZtau_2(a)(\omega^*_s\omega_s-1)\in\KOquad\label{KOunitary}
\end{gather}
for all $s \in [0,1]$. We set $z=u \oplus 1_3$, $\ee=\CALrho_2(1)$ and $\eep=1_6-\ee$, and abbreviate
\[
 C=\MZtaucos\ee,
\]
noting that this is a \cstar-algebra since $e$ and $\MZtaucos$
commute.
Note that $\MZtaucos =C + \eep\QQquad\eep\cong C+\QQsix$, and 
define $w=e\dot{z}e$ and
$w^\perp=\eep\dot{z}\eep$. We have that $[\ee,\dot{\omega}_s]=0$ by
\eqref{KOcommute}, so if we let
\[
\vv_s=\ee\dot{\omega}_s\ee\qquad \vvp_s=\eep\dot{\omega}_s\eep,
\]
we will get that
$\dot{\omega}_s=\vv_s\oplus\vvp_s$ for all $s$. Furthermore, $\vv_s$ is
a continuous path of unitaries by \eqref{KOunitary}, going from $w$ to
$\ee$ in $C$.
This implies that $\kone{w}{C}=0$, and consequently that 
$\kone{w}{Q}=0$. Since
$\kone{u}{\MMdo}=0$, we get
\[
0=\kone{\dot{z}}{Q}=\kone{w}{Q}+\kone{w^\perp}{Q}=
\kone{w^\perp}{Q}.
\]
We can now choose $n$ such that there exist a homotopy
\[
\xymatrix{{w^\perp \oplus n\eep=\vvp_0\oplus n\eep}\ar@{~}[rr]^-{w_s}&&{(n+1)\eep}}
\]
\SCOMMENT{mul?}in $\UG_{n+1}(\eep\QQquad\eep)$. Adding up ($v_s$, $\omega_s$, $n\ee$),  and identifying appropriately,
we
get a homotopy
\[
\xymatrix{{\dot{z}\oplus n 1_\QQquad}\ar@{~}[rr]^-{\Omega_s}&&{(n+1)1_\QQquad}}
\]
in $\UG_{n+1}( C+\eep\QQquad\eep)$. But since 
$ C+\eep\QQquad\eep =\MZtaucos$, this shows that
$\kone{\pi_6(u \oplus 1_3) }{\MZtaucos}=\kone{\dot{z} }{\MZtaucos}=0$.
Note that 
$D_{\MZtau_2}\cong \matrM_2(D_\MZtau)$ and
$\MZtaucos \cong \matrM_2 (\MZtauco)$, so that $\kone{\dot{u}}{\MZtauco}=0$.
 Consequently, applying the K-theory exact
sequence arising from \eqref{Drhotwo}, we may write
$
\kone{u}{K_1(\Drho)}=j_*(x)
$
for some $x\in K_1(\KOdo)$. 
\end{demo}
\begin{notes}
We are in fact using the Skandalis (nonunital) version of Paschke-Valette duality,
(cf.\ \cite{wlp:kcca}, \cite{av:rkg},
 \cite{gs:nnk},  \cite{nh:cetd}), in which the 
isomorphism
\[
\mu_1:\KK(A,B)\rightarrow K_1({\CALrho(A)\COM}\ee)
\]
is seen to map elements of the form $[\MZtau,\MZtau,u]$ to $[w]$, with
$\ee$ and $w$ as in
our proof above. Applying only the fact that this map is well-defined,
we conclude in our setting that $\kone{w}{\CALrho(A)\COM\ee}=0$.
\end{notes}

\subsubsection*{Obtaining inner automorphisms}
Recall that $\Mtau$ is an admissible scalar representation and that
$\MZtau$ denotes $0_2\oplus\Mtau:A\longrightarrow\MMdo$. We define,
for any $m\in\NN$, a 
\cstar-algebra $E_m$  by
\[
E_m=\{\mul{m}\MZtau(a)\mid a\in A\}+\KKK(H_{3m})\otimes
B+{\mathbb C}1_{3m},
\]
considered as a subalgebra of $M(\KKK(H_{3m})\otimes
B)$.

We need also to consider, for a unital \cstar-algebra $D$, the
group $\Aut_0(D)$ of automorphisms of $D$ which are connected to the
identity by a norm continuous path. 
It follows from \cite[8.6.12,8.7.8]{gkp:cag} that whenever $D$ is a unital
separable \cstar-algebra, $\Aut_0(D)\subseteq
\overline{\Inn(D)}$.

\begin{propo}{vII}
Let $A,B,\phi,\psi$ and $\Mtau$ be as in Proposition \ref{vI}. Whenever a finite subset
$\fset\subseteq A$ and $\epsilon>0$ is
given, there exist
$n\in\NN$ and a unitary $V\in
E_{n+1}$ such that 
\[
\norm{
V\left(
\left[\begin{matrix}\phit(a)&&\\&\Zha&\\&&\Mtau(a)\end{matrix}\right]\oplus
\mul{n}\MZtau(a)\right)V^*
-
\left[\begin{matrix}\psit(a)&&\\&\Zha&\\&&\Mtau(a)\end{matrix}\right]\oplus
\mul{n}\MZtau(a)}<\epsilon
\]
for all $a\in \fset$.
\end{propo}
\begin{demo}
Apply Proposition \ref{vI} to get $\Mut$, and choose $t_0$ with
$|1-t_0|<\epsilon/3$ and such that 
\begin{equation} \label{110}
\norm{
\Mut\left[\begin{smallmatrix}\phit(a)&&\\&\Zha&\\&&\Mtau(a)\end{smallmatrix}
\right]\Mut^*-
\left[\begin{smallmatrix}\psit(a)&&\\&\Zha&\\&&\Mtau(a)\end{smallmatrix}
\right]}<\epsilon/3
\end{equation}
for all  $a\in \fset$ and all $t>t_0$. Fix $t\in(t_0,1)$.
Because of Proposition \ref{trivkk}, Proposition
\ref{inKoneimage}
applies to $[\MZtau,\MZtau,\Mut]$. We can then choose a unitary
\[
w\in
\UG\left(\left[\begin{smallmatrix}\Zha&&\\&\KOo&\\&&\Zha\end{smallmatrix}\right]+\Odo\right)
\]
such that $\kone{\Mut}{K_1(\MDrho)}=\kone{w^*}{K_1(\MDrho)}$. Consequently,
there exists $n\geq 1$ and a homotopy
\[
\xymatrix{{\Mut \oplus w \oplus
1_{3n-3}}\ar@{~}[rr]^-{\omega_s}&&{1_{3n+3}.}}
\]
in $\UG_{n+1}(\MDrho)$. By our definition of $\MDrho$,
$\omega_s E_{n+1}\omega_s^*=E_{n+1}$, so we can define
a norm continuous family of automorphisms of $E_{n+1}$ by
$
\alpha_s=\Ad_{\omega_s}.
$
Clearly $\alpha_1=\id$, and since $\Ad_w\circ\MZtau=\MZtau$ because of
the special form of $w$, we get that $\alpha_0$ acts as
$\Ad_\Mut\oplus \mul{n}\id$ on elements of the form
\[
\left[\begin{smallmatrix}b&&\\&\Zha&\\&&\Mtau(a)\end{smallmatrix}\right]\oplus
\mul{n}\MZtau(a).
\]
Hence we have found
$\alpha\in\Aut_0(E_{n+1})$ such that
\begin{equation} \label{111}
\alpha\left(\left[\begin{smallmatrix}x&&\\&\Zha&\\&&\Mtau(a)\end{smallmatrix}\right]\oplus
\mul{n}\MZtau(a)\right)=\Mut\left[\begin{smallmatrix}x&&\\&\Zha&\\&&\Mtau(a)\end{smallmatrix}\right]\Mut^*\oplus
\mul{n}\MZtau(a)
\end{equation}
for all $a\in A$ and all $x\in \KO$.
As noted above, classical results about automorphisms allow us to find
$V \in E_{n+1}$ such that
$
\norm{
Va'V^*-
\alpha(a')}<\epsilon/3
$
for all  $a'\in \fset'$ given by
\[
\fset'=\left.\left\{\left[
\begin{smallmatrix}\phit(a)&&\\&\Zha&\\&&\Mtau(a)\end{smallmatrix}
\right]\oplus\mul{(n-1)}\MZtau(a)\right|a\in\fset\right\} .
\]
Since $\phit(a)\in\KO$, we get the desired
estimate from
\eqref{110} and \eqref{111}.
\end{demo}

As promised above we are now going to abandon our practice of not
identifying \cstar-algebras living on isomorphic Hilbert
spaces. Except in the proof of the following result, we will hence
drop the indices $H_1,H_2,\dots$, and talk only of the separable 
infinite Hilbert space $H$. 
We are going to consider 
\[
\Egen=\left.\left\{\left[\begin{smallmatrix}0&&\\&0&\\&&\genabsrep(a)\end{smallmatrix}\right]\right|
a\in A\right\}+\matrM_3(\KKK(H)\otimes 
B)+{\mathbb C}\left[\begin{smallmatrix}1&&\\&1&\\&&1\end{smallmatrix}\right] \subset \MMdo,
\]
where $\genabsrep$ is an absorbing (unital) representation of $A$ on
$\MMg$.
\begin{theor}{vIVPROTO}
Let $A$ be a unital, separable, nuclear C*-algebra and let $B$ be a  $\sigma$-unital \cstar-algebra.
Assume that $\phi$ and $\psi$ are two \hoos\ from
 $A$ to ${\mathcal K}(H) \otimes B$ which satisfy
$[\phi]_{KK}=[\psi]_{KK}$ in $KK(A,B)$.
 Let 
 $\genabsrep:A\longrightarrow \MMg$ be a unital
absorbing representation.
Whenever a finite subset
$\fset\subseteq A$ and $\epsilon>0$ is
given,  there  exists
a unitary $V\in \Egen$ such that 
\[
\norm{
V\left[\begin{matrix}\phit(a)&&\\&0&\\&&\genabsrep(a)\end{matrix}\right]V^*
-
\left[\begin{matrix}\psit(a)&&\\&0&\\&&\genabsrep(a)\end{matrix}\right]}<\epsilon
\]
for all $a\in \fset$.
\end{theor}
\begin{demo}
Apply Proposition \ref{vII} with some admissible scalar representation
$\Mtau$. By the fact that both $\mul{(n+1)}\Mtau$ and $\genabsrep$ are absorbing as
seen above, we can choose a sequence of unitaries
$u_m\in\LLL_B(H_{3n+2}\otimes B, H_{2}\otimes B)$ implementing the 
equivalence
$
0_1\oplus\Mtau\oplus\smul{n}{\MZtau}\abequiv 0_1\oplus\genabsrep.
$
Write $W_m=(1\oplus u_m)V(1\oplus u_m^*)$. With $m$ sufficiently
large, we  have
\[
\norm{W_m\left[\begin{smallmatrix}\phit(a)&&\\&\Zha&\\&&\genabsrep(a)\end{smallmatrix}\right]W_m^*
-
\left[\begin{smallmatrix}\psit(a)&&\\&\Zha&\\&&\genabsrep(a)\end{smallmatrix}\right]}<\epsilon,
\]
and since, by (ii) of Definition \ref{abequivalent}, $\Ad_{1\oplus
u_m}(E_{n+1})=\Egen$, we have that $W_m\in \Egen$. This proves the claim.
\end{demo}
\begin{notes}
The idea of proving results about
the automorphisms of the  \cstar-algebra $\Egen$  in order to
 obtain uniqueness results for morphisms originates with  Lin
(\cite{hl:saueh}) in the case $\genabsrep=\lind$.
\end{notes}

\subsection{Stable uniqueness}\label{doctorPROTO}
The importance of the auxiliary \cstar-algebra $\Egen$
used in Theorem \ref{vIVPROTO} above becomes apparent when one attempts to
employ our quasidiagonality condition \ref{quasidiagonalrep} to truncate the absorbing extensions in play to something more manageable
in terms of classification.
\subsubsection*{Achieving stable uniqueness}
We look at 
$\phi$ and $\psi$, two \hoos\ between unital \cstar-algebras $A$ to $B$.
To fix notation,
let $0$ denote the zero operator in $\MMg$ and let
$\genabsrep:A\longrightarrow \MMg$ be a quasidiagonal representation.
We consider a quasidiagonalization $(\genabsrep_n):A\longrightarrow
\matrM_{r_n}(B)$ by $(e_n)$, where we may and shall assume that  $(e_n)$
has the additional property that
\[
e_n\phi(a)e_n=\phi(a)\quad
e_n\psi(a)e_n=\psi(a). \quad n \geq 1.
\]

\begin{theor}{vVaumPROTO}
Let  $A$ be a unital, separable, nuclear \cstar-algebra and
let $B$ be a unital \cstar-algebra.
Assume that
there exists a quasidiagonal unital absorbing representation
$\genabsrep:A\longrightarrow \MMg$,
and let $(\genabsrep_n):A\longrightarrow \matrM_{r_n}(B)$
be a quasidiagonalization of  $\genabsrep$ by $(e_n)$ as above.

 Suppose that  $\phi,\psi:A\longrightarrow B$ are two 
$*$-homomorphisms with $\KKhom{\phi}=\KKhom{\psi}$ in $\KK(A,B)$,
such that $\phi(1)$ is unitarily equivalent to $\psi(1)$.
 Then for any
finite subset $\fset\subseteq A$ and any $\epsilon>0$, there exist
an integer $n$ and a unitary
$u\in\UG_{r_n+1}(B)$ satisfying
\[
\norm{
u\left[
\begin{matrix}
\phi(a)&\\
&\genabsrep_n(a)\end{matrix}\right]u^*-
\left[
\begin{matrix}
\psi(a)&\\
&\genabsrep_n(a)\end{matrix}\right]}<\epsilon
\]
for all $a\in\fset$.
Moreover we may arrange that $u(\phi(1)\oplus \gamma_n(1))u^*=
\psi(1)\oplus \gamma_n(1)$.
\end{theor}
\begin{demo}
After conjugating $\psi$ by a unitary in $B$ we may assume that
$\phi(1)=\psi(1)$.
We are going to compress by 
$
e_n'=e_n\oplus e_n\oplus e_n
$
(which is a quasicentral sequence in $\Egen$), in the conclusion of Theorem \ref{vIVPROTO}.
We have
\[
e_n'\left[\begin{smallmatrix}\phit(a)&&\\&0&\\&&\genabsrep(a)\end{smallmatrix}\right]
e_n'=
\left[
\begin{smallmatrix}
\phit(a)&&\\
&0&\\
&&\genabsrep_{n}(a)\end{smallmatrix}\right]
\]
and a similar equation for $\psit$. 
It is crucial to our argument that the unitary $V$
provided by Theorem \ref{vIVPROTO} satisfies  $\norm{[V,e_n']}\rightarrow
0$ because $V\in \Egen$. Therefore  by perturbing $e_n'Ve_n'$
to a unitary $v$ within $\UG_{3r_n}(B)$, for some large $n$,
we obtain:
\[
\norm{
v\left[
\begin{smallmatrix}
\phit(a)&&\\
&0&\\
&& \genabsrep_n(a)\end{smallmatrix}\right]v^*-
\left[
\begin{smallmatrix}
\psit(a)&&\\
&0&\\
&&\genabsrep_{n}(a)\end{smallmatrix}\right]}<\epsilon
\]
for all $a\in\fset$.
 Consider the projection
$e=\phi(1) \oplus \gamma_n(1)=\psi(1) \oplus  \gamma_n(1)$.
After a small perturbation of $v$ we may assume that $vev^*=e$.
 Then $w=eve$ is partial isometry
  in $\matrM_{r_n+1}(B)$ with $w^*w=ww^*=e$,
 and the unitary  
 $u=w+1_{r_n+1}-e \in \UG_{r_n+1}(B)$ will satisfy the conclusion of the Theorem.
\end{demo}

\section{Improved uniqueness results}\label{doctor}
This section contain successive refinements of uniqueness
results derived from the last result of the previous section.

\subsection{Stable uniqueness with bounds}

If one specializes to the case $\gamma=\lind$ in Theorem
\ref{vVaumPROTO}, one obtains a stable approximate unitary equivalence
of the form
\[
\norm{
u\left[\begin{matrix}\phi(a)&\\&\mul{n}\iota(a)\end{matrix}\right]u^*
-\left[\begin{matrix}\psi(a)&\\&\mul{n}\iota(a)\end{matrix}\right]
}<\epsilon.
\]
To make such a result useful in our quest to
classify \cstar-algebras, we need to refine our uniqueness results to the effect of
controlling the number $n$. More specifically, we need to know that these integers can be chosen
uniformly with respect to the targets; i.e.\ only depending on the
source algebra and, of course, the requirements on how closely the two
morphisms are to agree after composition by the unitary.

We also need to strengthen the theorem to allow for maps which are not
$*$-homomorphisms and only induce the same element locally in
$\Hom_\Lambda(\hugePK{A},\hugePK{B})$, rather than in $\KK(A,B)$.
To achieve such results, we are going to work with products of
\cstar-algebras, and we are going to depend on the results in Appendix
\ref{KKinj} regarding their $K$-theory.
If $A$ satisfies the UCT, then if follows from \eqref{UMCT} that
$ \Hom_\Lambda(\hugePK{A},\hugePK{B})$
 is isomorphic to  R{\o}rdam's group $KL(A,B)$ 
  \cite{mr:ccisa}.

\subsubsection*{Bounded stable uniqueness for $*$-homomorphisms}
\begin{theor}{vVIIum} 
Let $A$ be a simple, unital, nuclear, separable
 \cstar-algebra satisfying the UCT. 
Then for any finite subset $\fset\subseteq 
A$ and any $\epsilon>0$, there exists $n\in\NN$ with the
following property. For any admissible target $B$, any 
unital embedding $\iota:A\longrightarrow B$ and any 
pair of $*$-homomorphisms $\phi,\psi:A\longrightarrow B$ such that
${\phi}_*={\psi}_*$ in $\Hom_\Lambda(\hugePK{A},\hugePK{B})$, 
and $\phi(1)$ is unitarily equivalent to $\psi(1)$,
there exists a unitary $u\in\UG_{n+1}(B)$ such that
\[
\norm{
u\left[\begin{matrix}\phi(a)&\\&\mul{n}\iota(a)\end{matrix}\right]u^*
-\left[\begin{matrix}\psi(a)&\\&\mul{n}\iota(a)\end{matrix}\right]
}<\epsilon.
\]
for all $a\in \fset$.
Moreover we may arrange that $u(\phi(1)\oplus \mul{n}1)u^*=
\psi(1)\oplus \mul{n}1$.
\end{theor}
\begin{demo}
Suppose not and fix $\fset$ and $\epsilon$ for which the theorem
fails. Then for any $i$ we choose an admissible target algebra $B_i$ equipped with an
embedding $\iota_i:A\longrightarrow B_i$, and $\phi_i,\psi_i$ $*$-homomorphisms
with ${\phi_i}_*={\psi_i}_*$, and $\phi_i(1)$ unitarily equivalent
to $\psi_i(1)$, yet
\[
\inf_{u\in\UG_{i+1}(B)}\max_{a\in\fset}
\norm{
u\left[\begin{matrix}\phi_i(a)&\\&\mul{i}\iota_i(a)\end{matrix}\right]u^*
-\left[\begin{matrix}\psi_i(a)&\\&\mul{i}\iota_i(a)\end{matrix}\right]
}\geq \epsilon.
\]
We define $
\Phi,\Psi,I:A\longrightarrow\myPr
$ in the obvious way, and compose
with the canonical map to get
$\dotPhi,\dotPsi,\dotI:A\longrightarrow\myPrSu$. Since 
$\Phi$ and $\Psi$ induce the  same element of
$\prod\Hom_\Lambda(\hugePK{A},\hugePK{B_i})$ by construction, we get from
\PRODrefuhHOMinj\ that 
$\Phi_*=\Psi_*$ in
$\Hom_\Lambda(\hugePK{A},\hugePK{\myPr})$. Then of course also
$(\dotPhi)_*=(\dotPsi)_*$, and by \PRODrefuhPext\ we get that
$\KKhom{\dotPhi}=\KKhom{\dotPsi}$ in $\KK(A,\myPrSu)$.

Since
$\dotI:A\longrightarrow\myPrSu$ is a unital simple  embedding, we conclude
by Theorem \ref{vVaumPROTO} that there exist $n$ and a unitary
$w\in\UG_{n+1}(\myPrSu)$ intertwining $\dotPhi\oplus\mul{n}\dotI$ and
$\dotPsi\oplus\mul{n}\dotI$ up to $\epsilon$ on $\fset$. 
Let $u=(u_i) \in \UG_{n+1}(\myPr)$ be a unitary lifting $w$. Then
\[
\limsup_i \max_{ a \in \fset} \norm{
u_i\left[\begin{matrix}\phi_i(a)&\\&\mul{n}\iota_i(a)\end{matrix}\right]u_i^*
-\left[\begin{matrix}\psi_i(a)&\\&\mul{n}\iota_i(a)\end{matrix}\right]
}<\epsilon.
\]
yielding a contradiction after projecting onto $\matrM_{n+1}(B_i)$ for large $i$.
The last part of the proof is done exactly as the last part of the proof of 
 Theorem \ref{vVaumPROTO}.
\end{demo}
\begin{remar}{KL}
If $\phi$ and $\phi$ are as in the conclusion of either
Theorem \ref{vVaumPROTO} or
Theorem \ref{vVIIum}, it follows immediately from the definition of K-theory that
$\phi_*=\psi_*:\hugePK{A} \to \hugePK{B}$.
\end{remar}

\begin{remar}{simplifyI}
Under 
assumptions restricting the algebraic complexity on $\Ks(A)$ and
$\Ks(B)$ the result above can be simplified somewhat.
If we add, for instance, the assumptions that $K_0(A)$ be
torsion free and  $K_0(B)$ be
divisible, we need only require 
that $\phi_*=\psi_*$ on $\Ks(A)$. 
This is done basing the proof instead on injectivity of the maps
\begin{gather*}
{\Hom(\Ks({A}),\Ks({\mysmPr}))}\longrightarrow
{\prod\Hom(\Ks({A}),\Ks({B_i}))}
\\
{\KK(A,\mysmPrSu)}\longrightarrow \Hom(\Ks(A),\Ks(\mysmPrSu)).
\end{gather*}
We get the latter by
applying the UCT and note that since $K_0(A)$ is torsion free
\[
\Ext\left(K_0(A),K_1\left(\mysmPrSu\right)\right)=\Pext\left(K_0(A),K_1\left(\mysmPrSu\right)\right)=0
\]
from Corollary \ref{isac}(i), and that
\[
\Ext\left(K_1(A),K_0\left(\mysmPrSu\right)\right)=0
\]
since, along the lines of the first half of Lemma \ref{bdisalgco}, if
all $K_0(B_i)$ are 
divisible, then so is $K_0(\mysmPr)$.
\end{remar}

\subsubsection*{Bounded stable uniqueness for approximate morphisms}
We refer the reader to Appendix \ref{partialmaps}
for a discussion of partially defined maps on $\hugePK{-}$ and a
definition of $\huK$-triples.

\begin{theor}{vVII} 
Let $A$ be a simple unital, nuclear, separable
 \cstar-algebra satisfying the UCT.
 For any finite subset
$\fset\subseteq  
A$ and any $\epsilon>0$, there exist $n\in\NN$, and a $\huK$-triple $(\projsetKhu,\fseta,\delta)$  with the
following property.
For any admissible target $B$, and any three completely positive
contractions $\phi,\psi,\tau:A\longrightarrow B$ which are
$\delta$-multiplicative on $\fseta$, with $\tau$ unital and
$
\phi\pstar(p)=\psi\pstar(p)
$
in $\hugePK{B}$ for all $p\in\projsetKhu$, 
 and such that $\phi(1)$ and $\psi(1)$ are unitarily equivalent projections,
there exists a unitary
$u\in\UG_{n+1}(B)$ such that 
\[
\norm{
u\left[\begin{matrix}\phi(a)&\\&\mul{n}\tau(a)\end{matrix}\right]u^*
-\left[\begin{matrix}\psi(a)&\\&\mul{n}\tau(a)\end{matrix}\right]
}<\epsilon
\]
for all $a \in \fset$.
 One may arrange that
$u(\phi(1) \oplus \mul{n}1)u^*=\psi(1) \oplus \mul{n}1$.
\end{theor}

\begin{demo}
Seeking a contradiction we suppose that there is $\fset$ and
$\epsilon$ such that with $n=n(A,\fset,\epsilon)$ provided by Theorem
\ref{vVIIum}, no $\huK$-triple will
work. We choose sequences of $\huK$-triples $(\projsetKhu_i,\fseta_i,\delta_i)$ with the
properties 
\begin{rolist}
\item $\projsetKhu_i\subseteq \projsetKhu_{i+1}$,
and $\bigcup_{i\in\NN}\projsetKhu_i$ exhausts the semigroup 
\[
\bigcup_{m\in\NN}\bigproj(A\otimes C({\mathbb T})\otimes
C(W_m)\otimes\KKK)/\approx
\]
\item $\fseta_i\subseteq \fseta_{i+1}$,
$\overline{\bigcup_i{\fseta_i}}=A$.
\item $\delta_i>\delta_{i+1}$, $\delta_i\longrightarrow 0$.
\end{rolist}
By our assumption, we can then choose admissible targets $B_i$, which
we may assume are of the
same type, and
$\phi_i,\psi_i$ and $\tau_i$
which are $\delta_i$-multiplicative on $\fseta_i$ and satisfy
$
(\phi_i)\pstar(p)=(\psi_i)\pstar(p)
$
in $\hugePK{B}$ for all $p\in\projsetKhu$, and $\phi_i(1)$, $\psi_i(1)$ are unitarily equivalent
projections; yet
\[
\inf_{u\in\UG_{n+1}(B)}\max_{a\in\fset}\norm{
u\left[\begin{smallmatrix}\phi_i(a)&\\&\smul{n}\tau_i(a)\end{smallmatrix}\right]u^*
-\left[\begin{smallmatrix}\psi_i(a)&\\&\smul{n}\tau_i(a)\end{smallmatrix}\right]
}\geq \epsilon.
\]
Define
$
\Phi,\Psi,T:A\longrightarrow\myPr
$
from the sequences $(\phi_n),(\psi_n)$ and $(\tau_n)$, and compose
with the canonical map to get
$
\dotPhi,\dotPsi,\dotT:A\longrightarrow\myPrSu.
$
These maps are in fact $*$-homomorphisms by (ii) and (iii) above, so $\dotT$
provides a unital simple embedding of $A$ into $\myPr/\mySu$, and
$\dotPhi,\dotPsi$ induce
maps
\[
{\dotPhi}_*,{\dotPsi}_*:\hugePK{A}\longrightarrow\hugePK{\myPrSu}.
\]
We are going to show that ${\dotPhi}_*={\dotPsi}_*$.

We may check this on $p\in\projsetKhu_{j}$ by (i) above. Let 
$
C=
 C({\mathbb T})\otimes
C(W_m)\otimes\KKK,
$
with $m$ chosen appropriately. 
In the diagram
\[
\xymatrix{
0\ar[r]&
{K_0\left(\left(\mySu\right)\ampliC\right)}\ar[r]\ar[d]^{\cong}&
{K_0\left(\left(\myPr\right)\ampliC\right)}\ar[r]\ar[d]^{\eta}&
{K_0\left(\dfrac{\myPr}{\mySu}\ampliC\right)}\ar[r]\ar[d]^-{\dot{\eta}}&
0\\
0\ar[r]&
{\sum K_0\left(B_i\ampliC\right)}\ar[r]&
{\prod K_0\left(B_i\ampliC\right)}\ar[r]^{\kappa}&
{\dfrac{\prod K_0\left(B_i\ampliC\right)}{\sum K_0\left(B_i\ampliC\right)}}\ar[r]&
0}\]
$\eta$ is injective from \PRODrefuhETAinj\ because it can be naturally
identified with a component of $\ETA$. It is not hard
to check that the above diagram has exact rows. (The first row is induced by a quasidiagonal extension). By the five-lemma,  $\dot{\eta}$ is also injective, and hence it
suffices to show that
\[
\dot{\eta}[(\dot{\Phi}\otimes\id_C)(p)]=
\dot{\eta}[(\dot{\Psi}\otimes\id_C)(p)].
\]

With $\chi_0$ as in Appendix \ref{partialmaps},
lift $(\dot{\Phi}\otimes\id)(p)$ first to a  self-adjoint element
$({\Phi}\otimes\id)(p)$ in $\left(\myPr\right)\ampliC$, and then to a projection
$q$ in
\[
\chi_0(({\Phi}\otimes\id)(p))+\left(\mySu\right)\ampliC\subseteq \left(\myPr\right)\ampliC.
\]
Then $\eta([q])=([q_i])$ where $q_i=\chi_0((\phi_i\otimes\id)(p))$ for
all $i$ larger than some $i_{\Phi}$, since $\eta$ is induced by a family of
$*$-homomorphisms. We conclude that
\[
\dot{\eta}((\dot{\Phi}\otimes\id)_*(p))=
\kappa\eta([q])=
[\chi_0((\phi_i\otimes\id)(p))]_{i\geq i_{\Phi}}+\sum K_0(B_i\ampliC).
\]
Similarly, 
\[
\dot{\eta}((\dot{\Psi}\otimes\id)_*(p))=
[\chi_0((\psi_i\otimes\id)(p))]_{i\geq i_{\Psi}}+\sum K_0(B_i\ampliC),
\]
and these elements agree since the sequences coincide for $i\geq
i_\Phi,i_\Psi,j$.

Having proved that
 ${\dotPhi}_*={\dotPsi}_*$, and since $\dotPhi(1)$ is unitarily equivalent to 
$\dotPsi(1)$,
we may apply Theorem \ref{vVIIum} to find a unitary $w \in \UG_{n+1}( \myPrSu)$
such that
\[
\norm{
w\left[\begin{smallmatrix}\dotPhi(a)&\\&\smul{n}\dot{T}(a)\end{smallmatrix}\right]w^*
-\left[\begin{smallmatrix}\dotPsi(a)&\\&\smul{n}\dot{T}(a)\end{smallmatrix}\right]
}<\epsilon.
\]
for all $a\in \fset$.
Note that we may apply Theorem \ref{vVIIum} since
 whenever $B_i$ is a sequence of admissible
targets, $\myPrSu$ is an admissible target by \PRODrefuhCORO. 
We finish the argument by lifting the unitary as above.
\end{demo} 
\begin{remar}{simplifyII}
As in Remark \ref{simplifyI}, the premises of Theorem
\ref{vVII} simplify under extra 
assumptions on $\Ks(A)$ and $\Ks(B)$. If $K_0(A)$ is torsion free and
$K_0(B)$ is divisible, one needs only
check 
  that $\phi$ and $\psi$ agree on a $\Ks$-triple. This is because
divisibility passes from $K_0(B_i)$ to $K_0(\myPrSu)$ as outlined in Remark
\ref{simplifyI}.
\end{remar}
\section{Existence}\label{exist}
In this section we show how to realize locally a given $\KK$-element by a
difference of completely positive contractions (see Theorem \ref{EvI}). 
\subsection{Realizing group homomorphisms}
We refer the reader to Appendix \ref{partialmaps}
for the 
definition of $\huK$-triples and of the map $\chi_0$.

\begin{lemma}{clocom}
Let $E_i$ and $F_i$ be projections in $\MMg$ with\SCOMMENT{Want
$E_0-E_1$ small etc.}
\[
\norm{E_0-E_1}<\tfrac{1}{3}, \quad\norm{F_0-F_1}<\tfrac{1}{3}, \quad E_i-F_i\in\KOg
\]
for $i\in\{0,1\}$. Then $[E_0,F_0]=[E_1,F_1]$ in
$\KK_h(\CCC,B)\cong\KK(\CCC,B)$.
\end{lemma}
\begin{demo}
If $X_t=(1-t)E_0+tE_1$ and $Y_t=(1-t)F_0+tF_1$, then
$\spectrum(X_t),\spectrum(Y_t) \subseteq [0,1/3]\cup [2/3,1]$ for all $t\in[0,1]$.
Define $E_t=\chi_0(X_t)$ and $F_t=\chi_0(Y_t)$.
Then $(E_t,F_t)$, $0 \leq t \leq 1$ is a homotopy of Cuntz pairs
from $(E_0,F_0)$ to $(E_1,F_1)$. 
\end{demo}
\begin{lemma}{theEFconstant}
There exists $\lambda>0$ such that whenever $E,F$ are projections in
$\MMg$ with $E-F \in \KOg$
 (i.e.  $[E,F]$ defines a cycle in $\KK_h(\CCC,B)$), 
and $e$ is a projection in $\KOg$ with
\[
\norm{[e,E]}<\lambda\qquad\norm{[e,F]}<\lambda\qquad
\norm{e^\perp(E-F)e^\perp}<\lambda,
\]
then the natural isomorphism between $\KK_h(\CCC,B)$ and $K_0(B)$
takes $[E,F]$ to
\[
[\chi_0(eEe)]-[\chi_0(eFe)].
\]
\end{lemma}
\begin{demo}
Let $g=\chi_0(eEe),\, g'=\chi_0(e^\perp E e^\perp),\,
 h=\chi_0(eFe),\, h'=\chi_0(e^\perp F e^\perp)$.
For small $\lambda$, $g,g',h,h'$ are all projections and we may apply 
Lemma \ref{clocom} to get
\[ [E,F]=[g \oplus g',h \oplus h']=[g,h]+[g',h']=[g,h]+[g',g']=[g,h] \]
in $KK_h(A,B)$. Since $g,h \in \KOg$, the isomorphism
between $\KK_h(\CCC,B)$ and $K_0(B)$ takes $[g,h]$ to $[g]-[h]$.
\end{demo}

\begin{lemma}{makeitunital}
Let $A$ be a unital nuclear \cstar-algebra, and let $\epsilon>0$ and a
finite set $\fset\subseteq A$ be given. Then there is $\delta>0$ such
that for any unital \cstar-algebra $B$ and any completely positive
contraction $\phi:A\longrightarrow B$ which satisfies
$
\norm{\phi(1)^2-\phi(1)}<\delta,
$
then there exists a completely positive contraction $\psi$ with
$\norm{\phi(a)-\psi(a)}<\epsilon$ for all $a\in \fset$ such that $\phi(1)$
is a projection.
\end{lemma}
\begin{demo} 
Seeking a contradiction, we suppose that there exist $\fset$,
$\epsilon$ as well as sequences of unital \cstar-algebras $B_i$ and
completely positive contractions $\phi_i:A\longrightarrow B_i$ with
\[
\norm{\phi_i(1)^2-\phi_i(1)}\longrightarrow 0,
\]
yet for all completely positive contractions $\psi:A\longrightarrow
 B_i$ with $\psi(1)$ a projection, we have
\[
\sup_{a\in\fset}\norm{\phi_i(a)-\psi(a)}\geq \epsilon.
\]
If we let $\Phi=(\phi_i):A\longrightarrow\myPr$ and
 $\dot{\Phi}:A\longrightarrow\myPrSu$, then $\dot{\Phi}(1)=E$ is a
 projection. Let $(e_i)$ be a projection of $\myPr$ lifting $E$. This
 gives an isomorphism
\[
\gamma:\dfrac{\prod{e_iB_ie_i}}{\sum{e_iB_ie_i}}\longrightarrow
E\dfrac{\myPr}{\mySu}E,
\]
and by the unital version of the Choi-Effros theorem \cite{mdcege:cplpc} we can lift
$\gamma\inv\dot{\Phi}$ to a unital completely positive map
$\Psi=(\psi_i):A\longrightarrow \prod{e_iB_ie_i}$. Then
$\psi_i(1)=e_i$, and $\norm{\psi_i(a)-\phi(a)}$ tends to zero, leading to a contradiction.
\end{demo}

We now come to the main results in this section.
The reader is referred to Remark ~\ref{absqdnotes} for examples of quasidiagonal
absorbing representations.

Any element $\alpha \in \KK(A,B)$ induces a morphism $\alpha_*:
\hugePK{A} \to\hugePK{B}$.
If $p$ is a projection in $A$, and $\alpha$ is given by a
 Cuntz pair $(\tau,\gamma)$, then $\alpha_*$ takes
$[p] \in K_0(A) \cong \KK_h(\CCC, A)$ to $[\tau(p),\gamma(p)] \in
KK_h(\CCC,B)$. 

\begin{theor}{EvIPROTO}
Let  $A$ be a unital, separable, nuclear \cstar-algebra and 
let $B$ be a unital \cstar-algebra.
Assume that
there exists a quasidiagonal unital absorbing representation
$\genabsrep:A\longrightarrow \MMg$, 
and let $(\genabsrep_n):A\longrightarrow \matrM_{r_n}(B)$
be a quasidiagonalization of  $\genabsrep$ by $(e_n)$
as in Definition \ref{quasidiagonalrep}.

For any $\huK$-triple
$(\projsetKhu,\fset,\delta)$ 
there exist $N$ and a completely positive contraction
\[
\sigma:A\longrightarrow\matrM_{2r_N}(B)\qquad
\]
such that $\sigma$ and $\genabsrep_N$ are both $\delta$-multiplicative on $\fset$ and satisfy
\[
{\sigma}\pstar(p)-{(\genabsrep_N)}\pstar(p)=\alpha_*[p]
\]
for all $p\in\projsetKhu$.
\end{theor}
\begin{demo}
 For any admissible scalar representation 
$\Mtau$, if $\MZtau=0\oplus\Mtau$, then	
by \cite[2.6]{gs:nnk} $\alpha$ is represented by a \KK-cycle
$(\MZtau,\MZtau,x)$. By using the standard simplification given
by Proposition  17.4.3 in \cite{bb:koa} we may assume that $x$ is a  contraction.
Finally we may replace $x$ by the unitary 
\[ u =\left[\begin{smallmatrix}
 x & (1-xx^*)^{1/2} \\ -(1-x^*x)^{1/2} & x^*
\end{smallmatrix}\right] \]
Hence, in $\KK(A,B)$ we have
$ \alpha= [\MZtau, \MZtau,u]=[u\MZtau u^*,\MZtau,1]$.
This shows that for any admissible scalar representation 
$\Mtau$, $\alpha$ is represented by some Cuntz pair $(\rho, \MZtau): A \longrightarrow \MMg$.
We have 
\[ \rho(a)-\MZtau(a) \in \KOg \]
 for all $a \in A$.

Since  $\genabsrep$ is absorbing, we have $\Mtau\sim\genabsrep$ in the sense
of Definition \ref{abequivalent}. Dilating the unitaries trivially, we get a
sequence $u_i\in\UG(\MMg)$ with
\begin{equation}\label{fourteen}
\Zgenabsrep(a)-u_i\MZtau(a)u_i^*\in\KOg\qquad
\norm{\Zgenabsrep(a)-u_i\MZtau(a)u_i^*}\longrightarrow 0,
\end{equation}
where $\Zgenabsrep=0\oplus\genabsrep$. Note that for each $i$, $(u_i\rho u_i^*,u_i\MZtau u_i^*)$
is a Cuntz pair representing
$\alpha$. Then $X_i=(u_i\rho u_i^*,\Zgenabsrep)$ is also a Cuntz pair, and
$
[X_i]_*:\hugePK{A}\longrightarrow\hugePK{B}
$
converges pointwise to the group homomorphism induced by $\alpha$ as a
consequence of Lemma \ref{clocom}. Fix $i$ large enough that
$[X_i]_*[p]=\alpha_*[p]$ for each $p\in\projsetKhu$. 


Let $\tau= u_i\rho u_i^*$  and let $e_n \in \matrM_{r_n}(B)$
be as in the statement.
Then $f_n=e_n \oplus e_n \in matrM_{2r_n}(B)$ is an approximate unit of 
 projections  which in the obvious sense
quasidiagonalizes $\Zgenabsrep$ into $\gamma_n:A\longrightarrow \matrM_{2r_n}(B)$.
We set $\sigma_n(a)=f_n\tau(a)f_n$ and note that 
\[
[f_n,\tau(a)] \longrightarrow 0 \qquad
f_n^\perp (\tau(a)-\Zgenabsrep(a))f_n^\perp \longrightarrow 0
\]
for $a \in A$, since $\tau(a)-\Zgenabsrep(a) \in \KKK(H) \otimes B$.
We have that $ \projsetKhu\subseteq C$ where 
\[
C = \bigoplus_{m \leq M}{A\otimes C(\mathbb{T})\otimes C(W_m)\otimes\KKK}.
\]
Therefore for $N$  large enough, $\sigma_N$ is $\delta$-multiplicative
 on $\fset$ and 
\begin{gather*}
\norm{[f_N \otimes 1_C, (\tau \otimes \id_C)(p)]}<\lambda\qquad
\norm{[f_N \otimes 1_C, (\Zgenabsrep \otimes \id_C)(p)]}<\lambda\\
\norm{f_N^\perp \otimes 1_C( (\tau \otimes  \id_C)(p)- (\Zgenabsrep \otimes \id_C)(p))f_N^\perp \otimes 1_C}<\lambda,
\end{gather*}
for all $p \in \projsetKhu$.
By \ref{partialmaps} and Lemma \ref{theEFconstant} we have
\begin{eqnarray*}
\lefteqn{(\sigma_N )\pstar(p)- (\Zgenabsrep_N)\pstar(p)=}\\
& &[\chi_0(f_N \otimes 1_C (\tau \otimes \id_C)(p)f_N \otimes 1_C)]-
[\chi_0(f_N \otimes 1_C (\Zgenabsrep \otimes \id_C)(p)f_N \otimes 1_C)]=\\
& &[ (\tau \otimes \id_C)(p),(\Zgenabsrep \otimes
\id_C)(p)]_{KK_h}=[X_i]_*[p]= \alpha_*[p] 
\end{eqnarray*}
\end{demo}

Next we specialize the existence result    to the case
of a quasidiagonal source $A$. An application for purely infinite 
\cstar-algebras can be found 
in Theorem \ref{EvINFTY}.

A \cstar-algebra is RFD or \emph{residually finite-dimensional} if it
has a separating family of finite-dimensional representations. We say
that $A$ is \emph{locally RFD} if for any finite set $\fset$ and any
$\epsilon>0$, there exists an RFD subalgebra $A'$ of $A$, containing
all elements of $\fset$ up to $\epsilon$. 

\begin{theor}{EvI}
Let $A,B$ be  \cstar-algebras with $A$ nuclear, unital and quasidiagonal
and $B$ unital, and let $\alpha \in \KK(A,B)$.
For any $\huK$-triple
$(\projsetKhu,\fset,\delta)$ 
there exist $N$ and  completely positive contractions
\[
\sigma:A\longrightarrow\matrM_{N}(B)\qquad
\mu:A\longrightarrow\matrM_{N}(\CCC1_B)
\]
which are $\delta$-multiplicative on $\fset$ and satisfy
\[
{\sigma}\pstar(p)-{\mu}\pstar(p)=\alpha_*[p]
\]
for all $p\in\projsetKhu$.
We may arrange that $\sigma(1)$ and $\mu(1)$ are both
projections. Moreover if $A$ is locally RFD and if $\epsilon>0$
is given, we can arrange that
there is  
a unital RFD subalgebra $D$ of $A$ such that $\fset \subseteq_{\epsilon} D$
and the restriction of $\mu$ to $D$ is a $*$-homomorphism.
\end{theor}
\begin{demo}
Since $A$ is quasidiagonal, it has a quasidiagonal admissible
scalar representation $\Mtau$. Note that, using projections
$e_n\in\KKK (H)\otimes 1_B$ 
we obtain a quasidiagonalization consisting of maps
$\Mtau_n:A\longrightarrow\matrM_{r_n}(\CCC1_B)$.
The first claim follows from Theorem \ref{EvIPROTO}, by taking
$\mu=\Mtau_n$ for some large $n$.

Combining Lemma \ref{makeitunital} and Lemma \ref{close} we
may replace $\sigma$ and $\mu$ with completely positive
contractions which map the unit of $A$  to a projection. In case  $A$ is locally RFD, we find first an  unital RFD subalgebra $D$ of $A$ such that $F \subseteq_{\epsilon_1} D$. Then we work with an admissible
representation $\theta$ whose restriction to $D$ is a direct sum of
finite dimensional representations. It is then clear that one can choose
the quasidiagonalization such that the restriction of $\theta_n$ to $D$ is a $*$-homomorphism
for all $n$.
\end{demo}

\begin{notes}
An existence result for locally RFD \cstar-algebras has been
independently obtained by Lin \cite{hl:newtaf}. Our result is more
general and moreover it applies to the purely infinite \cstar-algebras. 
Nevertheless, the premises of the two results are the
same in the simple case,
 since it is proved in
\cite{bbek:gilfdc} that every simple, nuclear and quasidiagonal
\cstar-algebra is locally RFD,
\end{notes}
\section{Classification} 

In this section, we present applications of the uniqueness and
existence results to classification problems. The first part is
devoted to a class of finite \cstar-algebras which allow a further
refinement of the uniqueness result, leading to a complete
classification of \cstar-algebras in this class having $K_0$-group $\QQQ$. In
the last part, we show how the results also apply to reprove the 
classification theorem for purely infinite \cstar-algebras by
Kirchberg and Phillips.

\subsection{Approximate unitary equivalence}
In a class of \cstar-algebras studied by H.\ Lin it is
possible to absorb the stabilization required in Theorem \ref{vVII},
leading to further improved uniqueness results. Furthermore, this
class is contained in our class of admissible targets. In this section, we
develop these points.

\begin{defin}{deftaf} (\cite{hl:taf})
A simple unital \cstar-algebra is called TAF (tracially approximately finite dimensional) if for any finite subset $\fset \subseteq A$,  any $\ep >0$,
and any nonzero projection $q\in A$
there is a projection $p \in A$, $p\not= 1_A$,  and there is a finite dimensional
 \cstar-algebra\ $C \subseteq p^\perp A p^\perp$ with $p^\perp \in C$ \st\
\begin{rolist}
\item $\norm{[p,a]} < \ep $ for all $a \in \fset$.
\item $\dist(p^\perp a p^\perp, C)< \ep$ 
 for all $a \in \fset$.
\item $upu^*\leq q$ for some unitary $u\in A$.
\end{rolist}
\end{defin}

\begin{examp}{exoftaf}
 The real rank zero approximately (sub)homogeneous \cstar-algebras classified in
\cite{gaegg:ccrrzii} and \cite{mdgg:crahcrrz} are all TAF algebras.
Also the class of examples of nonnuclear subalgebras of AF algebras
 constructed in \cite{md:nsoaf} consists
entirely of TAF algebras.
\end{examp}
The conditions  (i)--(ii) from \cite{sp:lfdac} imply quasidiagonality of $A$. Let us
summarize a few other structural results on TAF \cstar-algebras that we
shall need. Lemmas \ref{tafweak} and  \ref{Lin} are from \cite{hl:taf}. 

\begin{lemma}{tafweak}(\cite[3]{hl:taf})
Let $A$ be a simple unital TAF \cstar-algebra. Then
$A$ has real rank zero,
stable rank one, and $K_0(A)$ is weakly unperforated in the sense of
\cite{bb:koa}. When  $p\in A$ is a
projection and $n\in\NN$, then both $pAp$ and $\matrM_n(A)$ are simple
unital TAF \cstar-algebras. 
\end{lemma}
\begin{lemma}{sublemma}
 Let $B$ be an infinite dimensional unital
separable simple C$^*$-algebra
of real rank zero, stable rank one, with $K_0(B)$ weakly unperforated.
Then for any $n \geq 1$ and any nonzero projection $f\in B$ there are mutually orthogonal
projections $e^1,\dots,e^n$ and $r$ in $B$ such that
$e^1+\dots+e^n+r=1_B$ with $[e^1]=\dots= [e^n]$ and $[r] < [f]$.
\end{lemma}
\begin{demo}
Let $QT(B)$ denote the normalized quasitraces on $B$.
The image of the natural map
$\rho_0:K_0(B)\longrightarrow \Aff(QT(B))$ is uniformly dense by
\cite[6.9.3]{bb:koa} and $K_0{B}$ has the strict ordering induced from 
$\rho_{0}$ by \cite[6.9.2]{bb:koa}. If $e$ is a projection we write $\hat 
e=\rho_{0}(e)$. By simplicity, we find $N$ big enough such that $N[f]>[1]$.
If $\ep=1/nN$, then $1/n-\ep >0$ and $n\ep \,\hat {1}<\hat {f}$.
Since the image of $\rho_{0}$ is uniformly dense, there is a 
projection $e \in B$ such that $(1/n-\ep)\hat {1} < \hat {e} < 1/n\, 
\hat {1}$. Therefore $0 < 1-n \hat{e} < n \ep \hat {1} < \hat {f}$, 
hence \[ 0 < [1]-n[e] < [f]. \]
This is readily seen to imply the statement.
Indeed if $d^i \in B$ are projections equivalent to $e$ 
then $d^1 \oplus \dots \oplus d^n$ is equivalent to a subprojection $d$ of $1_B$.
If $r=1-d$, then  $[r]=[1]-n[e]<[f]$.
\end{demo}

The TAF \cstar-algebras are prone to classification
because of the following factorization property.
\begin{lemma}{Lin} (\cite{hl:taf})
Let $A$ be a simple unital TAF \cstar-algebra.
For any $n \geq 1$, any finite subset $\fset\subseteq 
A$ and any $\epsilon>0$,
there are projections $p,q \in A$ with
 $p^\perp A p^\perp \cong \matrM_n(qAq)$ and $[p] \leq [q]$ and such that 
 there exists an approximate factorization
of $\id_A$
\[
\xymatrix{
A\ar[dr]_-{\nu}\ar@{=}[r]^-{\id_A}&A\\
&pAp\oplus \matrM_n(C)\ar@{^(->}[u]_-{\mu}}
\]
with $ \norm{\mu\nu(a)-a} < \epsilon$ for $a \in \fset$,
where $C$ is a unital finite dimensional \cstar-subalgebra of $qAq$,
$ \nu(a)=pap\oplus (\eta(a) \otimes 1_n)$
is $\epsilon$-multiplicative on $\fset$ with $\eta:A \longrightarrow C$ a unital
\CPAM, and $\mu$ is a unital $*$-monomorphism whose restriction to $pAp$ is the
natural inclusion.
\end{lemma}
\begin{demo}
We include a proof of this result which is somewhat different from
the original proof of Lin. We don't require $A$ to be nuclear.
We may assume that $A$ is infinite dimensional; it suffices to prove the statement 
with $8\epsilon$ instead of $\epsilon$.
Since $A$ is TAF we find a projection $P \in A$ and
a finite dimensional \cstar-algebra $C$ with $P^\perp \in C \subset P^\perp A P^\perp$ 
 such that for all $a \in \fset$ we have 
\begin{rolist}
\item $\|[P,a ]\| < \epsilon$ 
\item $P^\perp a P^\perp \in_{\ep} C$
\item
$(n+2)[P] \leq [1]$. 
\end{rolist}

The idea  of the proof is to find a unital embedding of $\matrM_{n}(\mathbb C) \oplus 
\mathbb C$ into the relative commutant of $C$ in $ P^\perp A P^\perp$
 such that the image of $\mathbb C$ is supported by a very 
small projection.
Write $C \cong \matrM_{m(1)} \oplus ... \oplus \matrM_{m(k)}$ and let $e_1,...,e_k$ be
the minimal central projections of $C$.
Let $B_i$ be the relative commutant
of $e_iCe_i \cong \matrM_{m(i)}$ in $e_iAe_i$. Then 
 $e_iAe_i \cong \matrM_{m(i)}(B_i)$ hence
$B_i$ is TAF being isomorphic to a corner of $A$. Let $f$ be a nonzero projection in $A$ with $(n+1)k[f] < [P]$. 
For each $1 \leq i \leq k$ we apply Lemma \ref{sublemma} for $B_i$. We obtain:
\[ e_i=e_i^1+\dots+e_i^n +r_i \]
where $e_i^j, 1 \leq j \leq n$ are mutually equivalent projections in $B_i$ and
$[r_i] < [f]$ in $K_0(A)$.
Set 
$e^j=e_1^j+\dots+e_k^j$, $1 \leq j \leq n$,
$r=r_1+\dots+r_k$ and
$e=e^1+\dots+e^n$.
Note that $P^\perp =e^1+\dots+e^n+r$
 with $e^j$ mutually equivalent
in the relative commutant of $C$ in $P^\perp A P^\perp$.
We have
 $(n+1)[r] =(n+1)([r_1]+\dots +[r_k]) <(n+1)k[f]< [P]$.
Therefore
\[ (n+1)([P]+[r]) < (n+2)[P] \leq [1] =[P]+n[e^1]+[r] \]
hence $n([P]+[r]) < n[e^1]$. By weak unperforation we get
$[P+r] \leq [e^1]$.
We are now ready to complete the proof. 
By Arveson's extension theorem, the inclusion $C \hookrightarrow A$ extends to
a completely positive contraction  $E: A \to C$. We have
$E(x)=x$ for $x \in C$ hence $\|a-E(a)\| \leq 2 \dist(a,C)$ for $a \in A$.  
Using (i) and (ii) we have for $a \in \fset$
\begin{eqnarray*}
 a &\sim_{2\ep}& PaP+P^\perp a P^\perp 
 \sim_{2\ep} PaP+E(P^\perp a P^\perp)\\&=& 
  PaP+E(P^\perp a P^\perp)r+ E(P^\perp a P^\perp)e 
  \sim_{4\ep} (P+r)a(P+r)+ E(P^\perp a P^\perp)e. 
\end{eqnarray*}
The last estimate  follows by compressing the estimate
$a\sim_{4\epsilon}PaP+E(P^\perp a P^\perp)$ by $P+r$.
It follows that $a \sim_{8\ep} 
   (P+r)a(P+r)+ E(P^\perp a P^\perp)e$.
We finish the proof by setting $p=P+r$, $q=e^1$ and noting that
 $E(P^\perp a P^\perp)e$
is of the form $\eta(a) \otimes 1_n$ since $e=e^1+\dots+e^n$ with
  $e^j$ mutually equivalent
in the relative commutant of $C$ in $P^\perp A P^\perp$.
\end{demo}

The uniqueness Theorems \ref{vVIIum} and \ref{vVII} apply to TAF algebras because of the following Proposition.
\begin{propo}{isTAFtarget}
A simple unital infinite-dimensional TAF \cstar-algebra is an admissible target algebra (of finite type).
\end{propo}
\begin{demo}
Let $B$ be a simple unital TAF \cstar-algebra. We get (i) and (ii) of Definition \ref{targetREST} by two results of
Rieffel (\cite{mr:dsrkc}, \cite{mar:hgugnt}). For \finiteref,
note that if $nx\geq 0$ then $nx+n[1_B]>0$ and $x+[1_B]>0$ by weak unperforation. Finally, to
prove \finitereftwo, assume that $\dim(B)=\infty$ and let $x$ and $n$
be given. The image of the natural map
$\rho_0:K_0(B)\longrightarrow\Aff(QT(B))$ is uniformly dense by
\cite[6.9.3]{bb:koa}, so we can find $z\in K_0(B)$ with
$\rho_0(x)-1<\rho_0(nz)<\rho_0(x)+1$. By \cite[6.9.2]{bb:koa} we have
$x-[1_B]\leq nz\leq x+[1_B]$ so $y=x-nz$ will work.
\end{demo}
\begin{theor}{unitaf}
Let $A$ be a simple unital, nuclear, separable
 TAF \cstar-algebra satisfying the UCT.
Then for any finite subset $\fset\subseteq 
A$ and any $\epsilon>0$, there exists a $\huK$-triple $(\projsetKhu,\fseta,\delta)$ with the
following property.
For any unital simple infinite-dimensional TAF \cstar-algebra $B$, and any two
unital completely positive contractions $\phi,\psi:A\longrightarrow B$ which are
$\delta$-multiplicative on $\fseta$, with
$\phi\pstar(p)=\psi\pstar(p)$ for all $p\in\projsetKhu$, there exists a unitary
$u\in \UG(B)$ such that 
\[
\norm{u\phi(a)u^*-\psi(a)}< \epsilon
\]
for all $a \in \fset$.
\end{theor}
\begin{demo}
Let us begin by outlining the proof. We first construct such a unitary
in the special case where $\phi,\psi$ are $*$-homomorphisms agreeing
on all of $\hugePK{A}$. This involves invoking Lin's factorization result \ref{Lin}
to pass to another pair of $*$-homomorphisms $\myphi,\mypsi$ which are
on a special form. Because Lin's result only gives an approximate
factorization, even though we start out
with $*$-homomorphisms, our proof will take us to a setting where our
uniqueness theorem for completely positive contractions Theorem \ref{vVII} is
needed. The general case will follow in the same way that
Theorem \ref{vVIIum} implies Theorem \ref{vVII}, by letting $n=0$ and
$T=0$ in the proof of Theorem \ref{vVII}. We include a sketch for the
benefit  of the suspicious reader.

{\sc Part 1:}  Given $\fset$ and
$\epsilon$, we are going to prove that
whenever
\begin{rolist}
\item $B=\myPrSu$
with each $B_i$ a unital simple infinite-dimensional TAF \cstar-algebra
(or $B$  itself a 
unital simple infinite-dimensional TAF \cstar-algebra).
\item $\phi,\psi:A\longrightarrow B$ are unital $*$-homomorphisms
\item $\phi_*=\psi_*: \hugePK{A} \longrightarrow \hugePK{B}$
\end{rolist}
then there exists $u\in\UG(B)$ with $\norm{u\phi(a)u^*-\psi(a)}<
\epsilon$ for all $a\in\fset$.

Let us thus fix 
$n$, $\projsetKhu$, $\fseta$ and $\delta$  by applying Theorem \ref{vVII}
to $\fset$ and $\epsilon/3$.
Furthermore, let $p$, $q$, $C$, $\nu$ and $\mu$ be given by Lemma \ref{Lin}
such that 
$\nu$ is $\delta$-multiplicative on $\fseta$ and 
$\norm{\mu\nu(a)-a} <\epsilon/3$ for all $a \in \fset$.

{\sc Step 1a:} Since $B$ has stable rank one by Lemma \ref{tafweak}, it 
has cancellation of projections, and because $\phi_*=\psi_*$ we may
assume,  after
conjugating $\psi$ by a unitary in $B$, that the 
restrictions
of $\phi$ and $\psi$ to the finite dimensional algebra $\mu(\CCC p \oplus \matrM_n(\CCC 1_C))$
are equal. Applying $\mu$ to the matrix units of $\matrM_n(\CCC 1_C)$ we can
define matrix units $(q_{ij})$ in $A$, where $q=q_{11}=1_C$. Let 
\begin{eqnarray}\label{efdef}
e=\phi(p)=\psi(p)\qquad f_{ij}=\phi(q_{ij})=\psi(q_{ij}),
\end{eqnarray}
abbreviating $f=f_{11}$. Invoking cancellation again, since $[p] \leq [q]$, we
can find a projection $e_0$ and a unitary $v$ in $B$ with
$v(e+e_0)v^*=f$. Let $g=e_0\oplus  1_B\in\matrM_2(B)$ and note that
$[g]=\mul{(n+1)}[f]$. Hence an isomorphism $\gamma:g\matrM_2(B)
g\longrightarrow \matrM_{n+1}(fB f)$ can be found. Denoting
the matrix units of $\matrM_{n+1}(\CCC f)$ by
$\tilde{f}_{ij}$ with $0\leq i,j\leq n $ we may choose $\gamma$ such that
\begin{eqnarray}\label{ongamma}
\gamma\left(\twobytwo{e_0b_0e_0}{}{}{ebe}\right)=\tilde{f}_{00}v(ebe+e_0b_0e_0)v^*\tilde{f}_{00}\qquad
\gamma\left(\twobytwo{0}{}{}{f_{ij}}\right)=\tilde{f}_{ij}\qquad
\end{eqnarray}
Combining all of this, we get a $*$-homomorphism $\myphi$ fitting in a diagram
\[
\xymatrix{ 
{A}\ar@{=}[r]^-{\id_A}\ar[dr]_-{\nu}&
{A}\ar[r]^-{\phi}&
{B} \ar@{^(->}[r]^-{\iota_2}&
{g\matrM_2(B)g}\ar[d]_-{\gamma}\\
&{pAp \oplus \matrM_n(C)}\ar@{^(->}[u]_-{\mu} \ar[rr]_-{\myphi}&&
{\matrM_{n+1}(fBf),}}
\]
where $\iota_2$ sends $B$ into the $(2,2)$ corner of
$g\matrM_2(B)g$. Identifying $\myphi$ using \eqref{efdef} and \eqref{ongamma} we get
that for all $d\in pAp$, $x\in\matrM_n(C)$,
\[
\myphi(d\oplus x)=\phi''(d)\oplus(\phi'\otimes \id_n)(x)
\]
where
\[
\phi':C\longrightarrow fB f\qquad\phi'':pAp\longrightarrow fBf
\]
are defined as corestrictions of $\phi$ and $\Ad_v\phi$, respectively. Furthermore, by symmetry of \eqref{efdef}, the
same procedure shows that
$\mypsi=\gamma\iota_2\psi\mu:pAp\oplus\matrM_n(C)\longrightarrow\matrM_{n+1}(fBf)$ has the form 
\[
\mypsi(d\oplus x)=\psi''(d)\oplus(\psi'\otimes \id_n)(x)
\]
for all $d\in pAp$, $x\in\matrM_n(C)$.

{\sc Step 1b:}  With $\iota_f:fBf\longrightarrow B$ we clearly have
that $(\iota_f\phi')_*=(\iota_f\psi')_*$ and
$(\iota_f\phi'')_*=(\iota_f\psi'')_*$. But since $f$ is full in $B$ as
the image of a full projection under a unital map, we get by
\cite{lgb:sihsc} that $\iota_f$ induces an isomorphism from
$\hugePK{fBf}$ to $\hugePK{B}$, so that $(\phi')_*=(\psi')_*$ and
$(\phi'')_*=(\psi'')_*$. From this we may assume, after conjugating
$\psi$ by a unitary of $fBf$, that $\phi'$ and $\psi'$ agree on $C$. Thus
the maps $\myphi\nu$ and $\mypsi\nu$ are of the form
\[
\myphi\nu(a)=\left[\begin{smallmatrix}
\phi''\omega(a)&&&\\&\phi'\eta(a)&&\\&&\ddots&\\&&&\phi'\eta(a)
\end{smallmatrix}\right]\qquad
\mypsi\nu(a)=\left[\begin{smallmatrix}
\psi''\omega(a)&&&\\&\phi'\eta(a)&&\\&&\ddots&\\&&&\phi'\eta(a)
\end{smallmatrix}\right]
\]
where we define $ \omega: A \longrightarrow pAp$ by $\omega(a)=pap$. Note that
 $\omega$ is $\delta$-multiplicative on $\fseta$. Now $\phi''\omega$,
 $\psi''\omega$  and $\phi'\eta$  
are $\delta$-multiplicative on $\fseta$ and by Lemma \ref{aha} we have
that 
 $(\phi''\omega)\pstar(p)=(\psi''\omega)\pstar(p)$ for all $p \in \projsetKhu$.
Therefore, Theorem \ref{vVII} applies
to the triple of maps $(\phi''\omega,\psi''\omega,\phi'\eta)$ if we can
prove that $fBf$ is an admissible target. When $B$ itself is a TAF
\cstar-algebra,
so is $fBf$ by Lemma \ref{tafweak}, and Proposition \ref{isTAFtarget} applies. When
$B=\myPrSu$  we note that there are projections $f_i \in B_i$ such that $fBf$ 
is isomorphic to
\[
\left.\prod{f_iB_i f_i}\right/\sum{f_iB_i f_i}
\]
and hence it is admissible of finite type by Lemma \ref{prosumofadmissible} since all the
$f_iB_i f_i$ are. 
Thus\SCHECK{}, by Theorem \ref{vVII}, there is a partial isometry $v \in  \matrM_{n+1}(fBf)$ such that
 $v^*v=\myphi\nu(1)$, $vv^*=\mypsi\nu(1)$ and 
 and $ \norm{v\myphi\nu(a)v^*-\mypsi\nu(a)}< \epsilon/3$ for all $a \in \fset$.
Since $\myphi\nu(1)=\mypsi\nu(1)=\gamma\iota_2(1)$ and 
$\gamma\iota_2(B)=\gamma\iota_2(1) \matrM_{n+1}(fBf)\gamma\iota_2(1)$,
we have that $v=\gamma\iota_2(u)$ for some unitary
$u \in B$.
Since $\gamma\iota_2$ is isometric 
we obtain
$ \norm{u \phi\mu\nu(a) u^*-\psi\mu\nu(a)}<\epsilon/3$, hence
$ \norm{u \phi(a) u^*-\psi(a)}<\epsilon$ for all $a \in \fset$,
 since $\norm{\mu\nu(a)-a} <\epsilon/3$ for all $a \in \fset$.

{\sc Part 2:} As the argument reducing to the case covered in {\sc
Part 1}
closely parallels that in the proof of Theorem \ref{vVII}, we only
sketch it here. If the theorem is false, we can choose sequences
$\projsetKhu_i,\delta_i,\fseta_i$ with the properties (i)--(iii) of that
proof, and corresponding simple unital TAF \cstar-algebras $B_i$ as
well as  unital completely positive contractions $\phi_i,\psi_i:A\longrightarrow B_i$
being $\delta_i$-multiplicative on $\fseta_i$, satisfying
$(\phi_i)\pstar(p)=(\psi_i)\pstar(p)$ for all $p\in\projsetKhu_i$; yet
\[
\inf_{u\in\UG(B)}\max_{a\in\fset}\norm{
u\phi_i(a)u^*
-\psi_i(a)
}\geq \epsilon.
\]
Define $\dot{\Phi},\dot{\Psi}:A\longrightarrow\myPrSu$ and check as in
the proof of Theorem \ref{vVII} that $\dot{\Phi}$ and $\dot{\Psi}$ are
$*$-homomorphisms inducing the same map on $\hugePK{A}$. Since
(i)--(iii) of {\sc Part 1} are met, we conclude from the first part of
the proof that there is a unitary $U\in\UG{(\myPrSu)}$  such that
$\norm{U\dot{\Phi}(a)U^*-\dot{\Psi}(a)}<\epsilon$ for all $a \in \fset$. Lifting
$U$ to a unitary $(u_i)\in\myPr$ and projecting, we get the
desired contradiction.
\end{demo}

\begin{notes}
As noted above,  Definition \ref{deftaf} is due to H.\ Lin.
It was motivated by a result of Popa \cite{sp:lfdac} and the classification theory
of $AH$ algebras (\cite{gaegg:ccrrzii},
\cite{mdgg:crahcrrz}). Definition \ref{deftaf} is in fact a version of the original definition, slightly simplified for
the class of simple algebras. To correlate this with
\cite{hl:taf} one compares to \cite[3.7]{hl:taf} and notes that since
Popa's conditions imply the (SP) property (all hereditary subalgebras
have a nonzero projection) it suffices to find a $p$ inside a generic corner
rather than inside a generic hereditary subalgebra.
\end{notes}

\subsection{Classification results}\label{tafclas}
We begin this section by presenting a shape type
classification result for simple TAF \cstar-algebras. Note that in
this setting, there is no need to appeal to our existence results, as
the \KK-classes are represented by completely positive
contractions from the outset. We then prove, this time using both
existence and uniqueness, that certain TAF \cstar-algebras are
isomorphic to the $AD$ algebra with the same $K$-theory. 
The main result in this section is Theorem \ref{classif}.
We refer the reader to 
\cite{acnh:dmakb} for the basics of asymptotic morphisms.
\begin{theor}{shape} Let $A$ and $B$ be two unital,
separable, nuclear, simple  TAF \cstar-algebras satisfying the UCT.
Suppose that
there are unital asymptotic morphisms $\phi=(\phi_t):A \longrightarrow B$ and  $\psi=(\psi_t):B \longrightarrow A$ such that 
$\phi_{*  }:\hugePK{A} \longrightarrow \hugePK{B}$ is bijective and $\phi_*^{\inv}=\psi_*$.
 Then $A$ is isomorphic to $B$.
\end{theor}
\begin{demo}
We may assume that both $A$ and $B$ are infinite-dimensional. We use repeatedly Theorem \ref{unitaf}
to find an increasing sequence of positive numbers $t_n \in (0,1)$ and sequences of unitaries $u_{2n-1} \in A$,
$u_{2n} \in B$, $n \geq 1$ such that the following diagram
is a two-sided approximate intertwining in the sense of Elliott \cite{gae:ccrrz}.
\[ \xymatrix{
A \ar[rr]^{\Ad u_1} \ar[dr]_{\phi_{t_1}}& 
& A \ar[rr]^{\Ad u_3} \ar[dr]_{\phi_{t_3}}& & A \ar[dr] \ar[r] & ... \\
& B \ar[rr]_{\Ad u_2} \ar[ur]^{\psi_{t_2}}& & B\ar[ur]^{\psi_{t_4}} \ar[rr]_{\Ad u_4}& & B \ar[r]& ...}
\]
\end{demo}

If $r$ is a positive integer we denote by 
$\DD{r}$ the \cstar-algebra of continuous functions
$f:[0,1] \longrightarrow \matrM_r(\mathbb C)$ such that  $f(0),f(1) \in \mathbb C 1_r$. 
We denote by $\DDn{r}$ the subalgebra of functions vanishing at
$0$. Let $\mathcal D$ be the class of algebras $B$  of the form
$ B=B_1\oplus\dots\oplus B_n$ where each $B_i$ is either a circle algebra 
$\matrM_k(C(S^1))$ or a dimension-drop algebra 
$\matrM_k(\DD{r})$.
An $AD$ algebra is a \cstar-algebra which is isomorphic to an inductive
limit of a sequence of \cstar-algebras in $\mathcal D$.

\begin{lemma}{inrealization}
Let $D$ be a dimension drop algebra $D=\mathbb I_{n}$ or $D=C_0(\RR)$.
Let $E$ be a simple unital TAF algebra and let $f$ be a nonzero projection in $E$.
Then the map $[D,fEf] \longrightarrow \KK(D,E)$ is surjective.
\end{lemma}
\begin{demo}
Let $n=1$ if $D=C_0(\RR)$. The proof uses the following series of facts:
\begin{rolist}
\item (\cite{mdtal:umctkg}) $\Dilim [D, \matrM_k(E)]=\KK(D,E)$. 
\item (\cite{mdtal:ecrrzc} when $D=\DDn{n}$) If $C$ is a finite dimensional
C*-algebra, and $\eta :D \longrightarrow C$ is a $*$-homomorphism, then
the map $ d \mapsto \eta(d) \otimes 1_n$ from $ D \longrightarrow \matrM_n(C)$ is
null homotopic.
\item  (\cite{tal:lsppc}) There is a finite subset $\fset \subseteq D$ and there is $\epsilon >0$ 
such that  if $\alpha,\beta:D \longrightarrow B$ are two $*$-\hos\ satisfying
$\norm{\alpha(d)-\beta(d)}<\ep$, then $\alpha$ is homotopic to $\beta$.
\item  (\cite{tal:lsppc}) For any 
 finite subset $\fset \subseteq D$ and $\epsilon >0$ 
 there is a finite subset $\fset_1 \subseteq D$ and there is $\epsilon_1 >0$ 
such that  if $\alpha:D \longrightarrow B$ is any \CPAM\ which is $\ep_1$-multiplicative 
on $\fset_1$, then there exists a $*$-\ho\ $\beta:D \longrightarrow B$ with
$\norm{\alpha(d)-\beta(d)}<\ep/2$ for all $d \in \fset$.
\end{rolist}

Fix $\fset$, $\ep$ as in (iii) and let $\fset_1$, $\ep_1$ be given by (iv). 
We may assume that $\fset \subseteq \fset_1$ and $ \ep > \ep_1$.
Let $x \in \KK(D,E)$. Then by (i) $x$ can be represented by some $*$-homomorphism $\gamma:D \longrightarrow  \matrM_k(E)$.
With $A= \matrM_k(E)$, $A$ is a simple unital TAF \cstar-algebra
by Lemma \ref{tafweak}. Since $E$ is simple there is $m \geq 1$ such that  $[1_E] \leq m[f]$.
Consider an approximation of $\id_{ A}$ 
provided by Lemma \ref{Lin} applied for the set $\gamma(\fset_1)$, $\ep_1/2$,
and the integer $mnk$:
\[
\xymatrix{
A\ar[dr]_-{\nu}\ar@{=}[r]^-{\id_A}&A\\
&pAp\oplus \matrM_{mnk}(C)\ar[u]_-{\mu}}
\]
Therefore \begin{equation} \label{twelwe}
 \norm{\mu\nu(a)-a} < \ep_1/2 
\end{equation}
for $a \in \gamma(\fset_1)$,
$C$ is a unital finite dimensional \cstar-subalgebra of $qAq$,
$ \nu(a)=\omega(a)\oplus (\eta(a) \otimes 1_{mnk})$, $\omega(a)=pap$, $\nu$
is $\epsilon_1$-multiplicative on $\gamma(\fset_1)$ and $\eta:A \longrightarrow C$ is a unital
\CPAM.
We may arrange that $\mu(p) \leq f$. 
Indeed, from $[p] \leq [q]$ and $[p]+mnk[q]=1$ we see that $(mnk+1)[p] \leq [1]$.
Since $\mu(1)=1_A$, we obtain $(mnk+1)[\mu(p)] \leq [1_A]=n[1_E] \leq mnk[f]$,
 hence 
$[\mu(p)] \leq n[\mu(p)] < [f]$ since $K_0(E)$ is weakly unperforated in the sense of \cite{bb:koa} by Lemma \ref{tafweak}.
After conjugating $\mu$ by a suitable unitary, we obtain $\mu(p) \leq f$.

Let us observe that $\omega\gamma$ and $\eta\gamma$ are 
 $\epsilon_1$-multiplicative on $\fset_1$.
By (iv) there are $*$-\hos\ $\omega':D \longrightarrow pAp$ and $\eta':D \longrightarrow C$ \st
if we set $\nu'=\omega' \oplus (\eta' \otimes 1_{mn})$, then
\begin{equation} \label {eleven}
\norm{\nu\gamma(d)-\nu'(d)}<\ep/2
\end{equation} for all $d \in \fset$.
{}From \eqref{twelwe} we have
 $\norm{\mu\nu\gamma(d)-\gamma(d)} < \ep_1/2 $
for $d \in \fset_1$. Combining this with \eqref{eleven} we get
\[ \norm{\mu\nu'(d)-\gamma(d)} < (\ep_1+\ep)/2<\ep  \]
for all $d \in \fset_1\cap \fset=\fset$.
By (iii) this implies that $\gamma$ is homotopic to $\mu\nu'$.
By (ii) $\nu'$ is homotopic to $\omega'$ so that $\gamma$ is homotopic to $\mu\omega'$. We conclude the proof by observing that
the image of $\mu\omega'$ is contained in $fEf$ since
$\mu\omega'(1)=\mu(p) \leq f$.
\end{demo}

\begin{lemma}{ifonethenK}
Let $A$ be a unital \cstar-algebra with $(K_0(A),[1_A])=(\QQQ,1)$.  Any
finite set of projections $\projsetnil\subseteq A\otimes\KKK$ can be
complemented to a $K_0$-triple $(\projset_0,\fseta,\delta)$ with the
property that for any unital completely positive
contraction $\phi:A \to A$ which is $\delta$-multiplicative on $\fseta$, one has
$\phi\pstar(p)=[p]$ for all $p\in\projset_0$. 
\end{lemma}
\begin{demo}
We may write $[p]=\tfrac rs [1_A]$, so
$\mul{s}p\oplus\mul{m}1_A\sim\mul{(r+m)}1_A$ for some $m \ge 0$. When $\phi$ is
sufficiently multiplicative, we have
\[
s\phi\pstar(p)+m[1_A]=\phi\pstar(\mul{s}p\oplus\mul{m}1_A)=
\phi\pstar(\mul{(r+m)}1_A)
=(r+m)[1_A].
\]
\end{demo}

The following theorem generalizes a result of Lin \cite{hl:taf}
by the fact that it allows non-zero (countable) $K_1$-groups.

\begin{theor}{reduction}
Let $A$ be a unital, separable, nuclear and simple TAF \cstar-algebra satisfying the UCT and
suppose that
$K_0(A) \cong \QQQ$ as ordered groups. Then $A$ is isomorphic to an
$AD$-algebra.
\end{theor}

\begin{demo}
Clearly $A$ is infinite-dimensional, and we may assume that $(K_0(A), [1_A]) \cong (\QQQ,1)$ as ordered pointed groups.
 For any finite subset $\fset \subseteq A$ and any $\ep>0$ we will find
an algebra $B \in \mathcal D$,
 and a $*$-\ho\ $\beta:B \longrightarrow A$ such that $\fset \subseteq_\ep \beta(B)$. 
This will prove the theorem as all elements of $\mathcal D$
are  semiprojective \cite{tal:lsppc}.
The class $\mathcal D$ was introduced before Lemma \ref{inrealization}.
As noted in Remark \ref{simplifyII}, applying Theorem \ref{vVII} 
to fixed $\fset$ and $\ep$ associated to 
\cstar-algebras 
$A$, $B$ with torsion-free divisible $K_0$-groups results
in a $\Ks$-triple $(\projset, \fseta,\delta)$ rather than a general
$\huK$-triple. 

Let $(\mathcal P_0, \mathcal G_0, \delta_0)$ and
 $(\mathcal V, \mathcal G_1, \delta_1)$ be a $K_0$-triple
 and a $K_1$-triple, respectively, given by Lemma
 \ref{star01} for the input $K_*$-triple 
 $(\mathcal P, \mathcal G, \delta)$. 
Let 
 $(\mathcal P', \mathcal G', \delta')$ be a $K_*$-triple
given by Lemma \ref{01star} for the input  triples
 $(\mathcal P, \mathcal G_0, \delta_0)$.
 and $(\mathcal V, \mathcal G_1, \delta_1)$.
We also may assume that 
 $ \mathcal G'$ and $ \delta'$ satisfy the conclusion of
 Lemma \ref{ifonethenK} 
applied for the unital \cstar-algebra $A$ and the set of projections $\projset_0$.

By \cite{gae:ccrrz} there is an $AD$ algebra $D$ and a group isomorphism
 $\kappa:K_1(A) \longrightarrow K_1(D)$.
By Theorem \ref{EvI} there exist \CPAMs\ $\sigma:A \longrightarrow M_N(D)$ and $\mu:A \longrightarrow
M_N(\CCC1_D)$ which are $\delta'$-multiplicative on $\fseta'$ and 
$\sigma\pstar(p')-\mu\pstar(p')=\kappa[p']$ for all $p' \in \fseta'$.
Here $\kappa$ is regarded as an element of
$\Hom(K_*(A),K_*(B)) \cong \Hom_\Lambda(\hugePK{A},\hugePK{B})$.
By the choice of the $K_*$-triple 
 $(\mathcal P', \mathcal G', \delta')$ 
 we have that
$\sigma\pstar(u)-\mu\pstar(u)=\kappa[u]$ for all unitaries $u \in \mathcal V$.
Note that $\mu\pstar(u)=0$ since $K_1(\CCC)=0$, so that we have
$\sigma\pstar(u)=\kappa[u]$ for all $u \in \mathcal V$.

Recall that $\sigma(1_A)=q$ is a 
projection, so that if we set $B=qM_N(D)q$, then $\sigma:A \longrightarrow B$
is a unital map.
Write $B$ as the inductive limit of an increasing sequence of algebras 
$B_k \in \mathcal D$, and let $j_k:B_k
\longrightarrow B$ be the inclusion map. Using the Choi-Effros
theorem as in Lemma 4.2 of \cite{mdtal:kasa} we find a sequence of \CPAMs\
 $\eta_k : B \longrightarrow B_k$ such that 
$j_k\eta_k $ converges to $\id_B$ in the point-norm topology.
Choose $k$ large enough so that 
\begin{equation} \label{AAA}
\sigma\pstar(u)=(j_k\eta_k\sigma)\pstar(u)
\end{equation}
for all $u \in \mathcal V$.

Consider the group morphism
$\kappa^{-1}(j_k)_*:K_1(B_k) \longrightarrow K_1(A)$.
There is a  a unital $*$-\ho\
$\gamma :B_k \longrightarrow A$ 
such that $\gamma_*=\kappa^{-1}(j_k)_*$,
obtained as follows. Write $B_k=\matrM_{\ell(1)}(B_{k1}')\oplus\dots\oplus 
\matrM_{\ell(r)}(B_{kr}')$ where the $B_{ki}'$ are either $C(S^1)$ or $\DD{n(i)}$.
Let $q_1+\cdots +q_r=1_A$ be a partition of $1_A$ by nonzero projections.
We have $K_0(q_iAq_i) \cong K_0(A) \cong \QQQ$ as ordered groups.
Therefore we find mutually orthogonal
projections $p_i$ such that $\ell(i)[p_i] \cong [q_i]$.
Since $A$ has cancellation of projections we find a unital inclusion
\[  \bigoplus_{i=1}^r M_{\ell(i)}(p_iAp_i) \subseteq A .\]
Using Lemma \ref{inrealization} we find unital $*$-\hos\
$\gamma_i :B_i \longrightarrow p_iAp_i$ such that 
$\gamma=\oplus_{i=1}^r (\gamma_i \otimes \id_{\ell(i)})$ has the desired property. 

Next we want to show that $\gamma \circ(\eta_k\sigma)$ gives an approximate
factorization of $\id_A$ on $K_*(A)$.
More precisely we want that
$(\gamma\eta_k\sigma)\pstar(p)=[p]$ for all $p \in \projset$.
By virtue of our choice of the $K_0$ and $K_1$-triples above, it suffices to
show that
$(\gamma\eta_k\sigma)\pstar(p_0)=[p_0]$ for all $p_0 \in \projset_0$ and
$(\gamma\eta_k\sigma)\pstar(u)=[u]$ for all $u \in \mathcal V$.
Using Lemma \ref{aha} twice, the definition of $\gamma$ and \eqref{AAA} we have
\[ 
(\gamma\eta_k\sigma)\pstar(u)=
\gamma_*(\eta_k\sigma)\pstar(u)=
\kappa^{-1}(j_k)_*(\eta_k\sigma)\pstar(u)=
\kappa^{-1}(j_k\eta_k\sigma)\pstar(u)=\kappa^{-1}\sigma\pstar(u)= [u] 
\]
for all $u \in \mathcal V$.
It remains to check that
$(\gamma\eta_k\sigma)\pstar(p_0)=[p_0]$ for all $p_0 \in
\projset_0$, but
this follows from Lemma \ref{ifonethenK} by our choice of
the $\mathcal G'$ and $\delta'$.

Define $\alpha=\eta_k\sigma : A \longrightarrow B_k$.
We have seen that $(\gamma\alpha)\pstar(p)=[p]$ for all $p \in \projset $.
Therefore by Theorem
 \ref{vVII} there is unitary $u \in A$ such that if $\beta=\Ad u \circ \gamma$,
then $\norm{\beta \alpha(a)-a}< \ep$ for all $a \in \fset$,
 hence $\fset \subseteq_\ep \beta(B)$. 
\end{demo}
\begin{theor}{classif}
Let $A, B$ be  unital, separable, nuclear
and simple  TAF \cstar-algebras satisfying the UCT. Suppose that
$(K_0(A), [1_A]) \cong (K_0(B), [1_B]) \cong (\QQQ, 1)$ as ordered groups,
and $K_1(A) \cong K_1(B)$. Then $A$ is isomorphic 
to $B$.
\end{theor}
\begin{demo}
By Theorem \ref{reduction} both $A$ and $B$ are isomorphic to simple
real rank zero $AD$ algebras. These are classified by their $K$-theory data
as proved by Elliott \cite{gae:ccrrz}.
\end{demo}
\begin{notes}
Lin proved Theorem \ref{reduction} for $K_1(A)=0$ in \cite{hl:taf}. He
subsequently, independently from and at about the same time as the
present work, generalized his result to allow general $K_1$-groups in
\cite{hl:newtaf}. This paper also discusses classification of
\cstar-algebras with other $K_0$-groups.
\end{notes}
\subsection{Purely infinite \cstar-algebras}\label{piclas}
The purpose of this section is to demonstrate how our methods can be
applied to give the classification result of Kirchberg and Phillips
starting from three basic, albeit deep, structural results about purely
infinite \cstar-algebras.
The methods used here are very similar to those used in the finite case,
with Cuntz' algebra $\OOO_2$ playing the role of 
$\matrM_n(\CCC)$. The exposition will emphasize this similarity.
 To make it very clear exactly how much we
need to import from the theory of this class of \cstar-algebras we
collect the required results below.

\begin{Rolist}
\item (\cite[2.8]{ekncp:eeccfco}) Any exact, separable and unital
\cstar-algebra embeds unitally into $\OOO_2$.
\item (\cite[3.6]{mr:cilca})
Let $B$ be an admissible target algebra. If
$\phi,\psi:\OOO_2\longrightarrow B$ is a pair of unital $*$-homomorphisms,
then for any finite set $\fset\subseteq \OOO_2$ and any $\epsilon>0$ there exists
$u\in\UG(B)$ with 
$ \norm{u\phi(a)u^*-\psi(a)}< \epsilon$ for all $a\in\fset$.
\item (A variation of \cite[2.4]{ncp:auehoca})
Let\SCHECK{} $A$ be a purely infinite nuclear separable
unital \cstar-algebra. For any $n \geq 1$, any finite subset $\fset\subseteq 
A$ and any $\epsilon>0$,
there is a projection $p \in A$ with
 $p^\perp A p^\perp \cong \matrM_n(\OOO_2)$ such that 
 there exists an approximate factorization
of $\id_A$
\[
\xymatrix{
A\ar[dr]_-{\nu}\ar@{=}[r]^-{\id_A}&A\\
&pAp\oplus \matrM_n(\OOO_2)\ar@{^(->}[u]_-{\mu}}
\]
with $ \norm{\mu\nu(a)-a} < \epsilon$ for $a \in \fset$,
$\nu(a)=pap\oplus (\eta(a) \otimes 1_n)$
is $\epsilon$-multiplicative on $\fset$ with $\eta:A \longrightarrow \OOO_2$ a unital
\CPAM, and $\mu$ is a unital $*$-monomorphism
 whose restriction to $pAp$ is the
natural inclusion.
\end{Rolist}

\begin{remar}{commentstoimport}
(I) is needed only for nuclear algebras.
We have rephrased (II) and (III) slightly to suit our needs. R\o rdam
requires that $B$ satisfies
\[
\operatorname{cel}(B)<\infty\qquad \UG(B)/\UG_0(B)\cong K_1(B)
\]
and this follows by Definition \ref{targetREST} as
explained in Appendix \ref{KKinj}. Also, Phillips proves (iii) only
for $n=1$, but if we  write $\mu:\OOO_2\longrightarrow 
\matrM_n(\OOO_2)$ for the map sending $x$ to $\diag(x,\dots,x)$ and
choose a unital embedding
$\iota:\matrM_n(\OOO_2)\longrightarrow \OOO_2$, then (II)
shows that a unitary $u\in\UG(\OOO_2)$ exists making the diagram
\[
\xymatrix{
{\OOO_2}\ar_-{\mu}[dr]\ar[rr]^-{\Ad u}&&{\OOO_2}\\
&{\matrM_n(\OOO_2)}\ar[ur]_-{\iota}}
\]
commute up to $\epsilon$ on $\fset$. Thus we may replace $\OOO_2$ by
$\matrM_n(\OOO_2)$ without loss of generality (and without using the
theorem we are aiming for).
\end{remar}

Apart from (I)--(III), all we need to know about a purely infinite
\cstar-algebra $A$ is that it  has
real rank zero (\cite{sz:ppisc}), that the canonical maps from $\bigproj(A)$ and $\UG(A)$
to $K_0(A)$ and $K_1(A)$ are surjective (\cite{jc:kcc}), and
furthermore, that if $p,q\in A$
are nonzero projections, then
\[
[p]=[q]\infer p\sim q.
\]
The latter fact is also from \cite{jc:kcc}.
Note that this shows:

\begin{propo}{isINFTYtarget}
A purely infinite simple unital \cstar-algebra is an admissible target (of
infinite type).
\end{propo}
\begin{theor}{EvINFTY}
 Let $A,B$ be  unital \cstar-algebras with $A$ nuclear separable   and
$B$ containing a unital copy of $\OOO_2$, and let $\alpha \in \KK(A,B)$.
Then for any $\huK$-triple
$(\projsetKhu,\fset,\delta)$ 
there exists a completely positive contraction
$ \sigma:A\longrightarrow B $
which is $\delta$-multiplicative on $\fset$ and satisfies
$
{\sigma}\pstar(p)=\alpha_*[p]
$
for all $p\in\projsetKhu$.
\end{theor}

\begin{demo}
 By
(I), $A$ embeds unitally into $\OOO_2$, so with
$\iota$ defined as the composite
\[
\xymatrix{{A}\ar@{^(->}[r]&{\OOO_2}\ar@{^(->}[r]&{B}}
\]
we obtain a unital simple embedding. By Lemma \ref{siabsorbs}, the
representation $\lind$ is absorbing, and it is clearly also
quasidiagonal as it commutes with the  projections $e_n=\mul{n}1_B$. By Theorem
\ref{EvIPROTO}, we find $N$ and 
$\sigma:A\longrightarrow \matrM_{2N}(B)$ such that 
\[
(\sigma )\pstar(p)- (\gamma_N)\pstar(p)=\alpha_*[p]
\]
where $\gamma_N$ is a
$*$-homomorphism of the form
\[
\gamma_N(a)= 0\oplus \mul{N}\iota(a),
\]
so that $(\gamma_N)_*=0$ by the fact that $\iota$ factors through
$\OOO_2$. Thus $(\sigma)\pstar=\alpha_*$ on $\projset$, and assuming, as we may (by Lemma \ref{makeitunital}), that $\sigma(1)$ is a projection $p$,
we have that $p\in\matrM_N(B)$ is a subprojection of $\mul{N}1_B$ which
is equivalent to a subprojection of $1_B=1_{\OOO_2}$ via some unitary $u$. We may hence
replace $\sigma$ by   $u\sigma u^*:A\longrightarrow B$ inducing the same
partial map on $\projsetKhu$.
\end{demo}
\begin{remar}{kirchstrong}
Kirchberg proved a stronger form of the previous theorem where $\alpha$ is lifted
to a $*$-monomorphism.
\end{remar}
\begin{theor}{uniINFTY}
Let $A$ be a  purely infinite separable unital nuclear \cstar-algebra satifying the UCT.
Then for any finite subset $\fset\subseteq 
A$ and any $\epsilon>0$, there exists a $\huK$-triple
$(\projset,\fseta,\delta)$  with the 
following property.
For any purely infinite simple unital \cstar-algebra $B$, and any two
unital completely positive contractions $\phi,\psi:A\longrightarrow B$ which are
$\delta$-multiplicative on $\fseta$, with
$\phi\pstar(p)=\psi\pstar(p)$ for all $p\in\projsetKhu$, there exists a unitary
$u \in \UG(B)$ such that $\norm{u\phi(a)u^*-\psi(a)}< \epsilon$
for all $a \in \fset$.
\end{theor}
\begin{demo}
The proof follows closely that of Theorem \ref{unitaf}, and we are only going to
indicate the changes needed. Here {\sc Part 1} of the proof deals with
a pair of $*$-homomorphisms into $B=\myPrSu$ with each $B_i$ a purely
infinite \cstar-algebra. Applying first Theorem \ref{vVII} and then
(III) we get $n,\projsetKhu,\fseta,\delta,p,q,\nu$
and $\mu$ as in that proof. Then {\sc Step 1a} applies verbatim as
soon as one notes that because they are all nonzero, one has cancellation on all the projections in
play by the result of Cuntz mentioned above. In {\sc Step 1b}
one applies (II) to obtain a unitary of $fBf$ which conjugates
$\psi'$ to $\phi'$ to within $\epsilon/3$ on
$\fset$. \SCOMMENT{Careful now!}This suffices to
achieve the desired conclusion by the argument given in Theorem \ref{vVII}. Furthermore, to get that $fBf$
is admissible of infinite type, one uses the same argument to
reduce to the case of showing that  a corner of a purely infinite
\cstar-algebra is also admissible of infinite type. This is clear
since it is itself purely infinite. Finally, {\sc
Part 2} of the proof carries through verbatim because of Proposition \ref{isINFTYtarget}. 
\end{demo}

\begin{theor}{pisunclas} (\cite{ek:cpicukt},
\cite{ncp:ctnpisc})
Let $A$ and $B$ be purely infinite separable unital nuclear
\cstar-algebras satisfying the UCT. 
Then any isomorphism
$\kappa:(\Ks(A), [1_A])\longrightarrow (\Ks(B), [1_B])$ is induced by a $*$-isomorphism.
\end{theor}
\begin{demo}
We may assume that $A$ and $B$ are in the standard form, that is  $1_A$ and $1_B$ both
represent the zero class in their respective $K_0$-groups.
It follows by \cite{jc:kcc} that $A$ and $B$ both contain unital
copies of $\OOO_2$. We may thus apply the existence result Theorem \ref{EvINFTY} to get
maps $\sigma_i:A\longrightarrow B$ and $\tau_i:B\longrightarrow  A$
which are increasingly multiplicative on larger and larger sets, and
induce $\kappa$ and $\kappa\inv$, respectively, on larger and larger
subsets of $\hugePK{A}$ and $\hugePK{B}$. Arranging this appropriately, we may
conclude by the uniqueness result Theorem \ref{uniINFTY} that
unitaries $u_i,v_i$ exist making
\[ \xymatrix{
A \ar[rr]^{\Ad u_1} \ar[dr]_{\sigma_{1}}& 
& A \ar[rr]^{\Ad u_2} \ar[dr]_{\sigma_{2}}& & A \ar[dr] \ar[r] &{\ldots} \\
& B \ar[rr]_{\Ad v_1} \ar[ur]^{\tau_{1}}& & B\ar[ur]^{\tau_{2}} \ar[rr]_{\Ad v_2}& & B \ar[r]& {\ldots}}
\]
an approximate intertwining in the sense of Elliott
(\cite{gae:ccrrz}). 
\end{demo}
\appendix
\section{Appendix on fine $K$-theoretical points}\label{KKap}
\subsection{$K$-theory of products}\label{KKinj}
It is well known, although perhaps not as well known as it should be,
that the natural map $\Ks(\myPr)\to\prod\Ks(B_i)$ is not an isomorphism in
general. In this section
we are going to study, for large classes of
\cstar-algebras $B_i$, injectivity properties of
the natural 
map 
\[
\ETA:\hugePK{\myPr}\longrightarrow \prod{\hugePK{B_i}}
\]
defined by collecting the maps induced by the projections
$\pi_i:\myPr\longrightarrow B_i$. We then use this information to prove that
surprisingly often, $\Pext(-,\Ks(\myPrSu))$ will vanish.
\subsubsection*{Five quantities}
We are going to define five quantities $\vali$ (``cancellation order''),
$\valii$ (``perforation order''), $\valiii$, $\valiv$ (``element
lifting order'') and
$\valv$ (``infinite height perturbation order'') in $\NN\cup\{\infty\}$  for any unital \cstar-algebra $B$ by
declaring
\begin{list}{}{\setlength{\leftmargin}{1.99cm}}
\item[$\vali(B)\leq\ell$] if whenever $p,q\in B\otimes\KKK$, then $[p]=[q]\Longrightarrow p\oplus \mul{\ell}1_B\sim q\oplus
\mul{\ell}1_B$
\item[$\valii(B)\leq \ell$] if for any $x \in K_0(B)$ such that $nx \geq 0$ for some $n>0$,
one has $x+\ell [1_B] \geq 0$.
\item[$\valiii(B)\leq \ell$] if the canonical map $\bigproj(\matrM_\ell(B))\longrightarrow K_0(B)$ is
surjective
\item[$\valiv(B)\leq \ell$] if the canonical map $\UG_{\ell}(B)\longrightarrow K_1(B)$ is surjective
\item[$\valv(B)\leq\ell$] if 
for any $x$ in $K_0(B)$  and any  $n\not=0$, there is
$y$ in $K_0(B)$ such that $-\ell[1_B] \leq y \leq  \ell[1_B]$ and $x-y\in nK_0(B)$.
\end{list}
and declaring the value to be $\infty$ when no such $\ell$ exists.
When $\{B_i\}_{i\in I}$ is a family of \cstar-algebras, we define 
$\vali(\{B_i\})=\sup_{i}\vali(B_i)$,
 and so forth.
We also write $\rera(\{B_i\})=0$ when each $B_i$ has real rank zero.

The following results -- some of which are known, cf.\
\cite{setal:ccsr} -- demonstrate the relevance of these quantities to
the map $\ETA$. We denote the unit of $B_i$ by $1_i$ and the units of $\myPr$ and
$\myPrSu$ by $\unPr$ and $\unPrSu$, respectively. 

\begin{lemma}{aboutetanil}
Let  $B_i$ be a sequence of unital \cstar-algebras for which
$\vali(\{B_i\})<\infty$. Then
$\eta^0$ is injective, and 
the image of $\eta^0$ equals
\[
\left.\left\{(x_i)\in \prod K_0(B_i)\right| \exists M\in\NN\forall i\in\NN:
-M[1_{i}]\leq x_i\leq {M}[1_{i}]\right\}.
\]
Furthermore, if $\eta^0(x)=(x_i)$ with each $x_i\geq 0$, then  $x\geq0$.
If in addition $\valii(\{B_i\})<\infty$, then $\Im\eta^0$ is a pure
subgroup of $\prod K_0(B_i)$. If $\valiii(\{B_i\})<\infty$, then $\eta^0$ is surjective.
\end{lemma}
\begin{demo}
Assume that $\vali(\{B_i\})\leq \ell$.
To prove injectivity, let $x=[(p_i)]-[(q_i)]\in K_0(\myPr)$ be given by
$p_i,q_i\in\matrM_N(B_i)$  and assume that  $\eta^0(x)=0$. We have
that $p_i\oplus \mul{\ell}1_{i}\sim
q_i\oplus \mul{\ell}1_{i}$, wherefrom
$(p_i)\oplus\mul{\ell}\unPr\sim
(q_i)\oplus\mul{\ell}\unPr$, proving $x=0$.
To prove that the image is contained in the set of bounded sequences, write $x\in K_0(\myPr)$ as
 $x=[(p_i)]-[(q_i)]$ with $p_i,q_i\in\matrM_N(B_i)$.
 By definition of positivity,
$-{N}[\unPr]\leq x\leq {N}[\unPr]$ and we
get the result by applying $\pi_i$. For the other
inclusion, assume that $x_i\in K_0(B_i)$ is given with ${-M} [1_{i}]\leq
x_i\leq {M}[1_{i}]$. We can write $x_i+{M}[1_{i}]=[r_i]$ for some
projection $r_i$ in $B_i\otimes\KKK$, and since $[r_i]\leq{2M}[1_{i}]$ we get
that $2M[1_{i}]=[r_i\oplus s_i]$ for some projection $s_i$ in  $B_i\otimes\KKK$.
By $\vali(\{B_i\})\leq\ell$ we see that 
$(2M+\ell)\cdot 1_i\sim r_i \oplus s_i\oplus \mul{\ell}1_i$,
hence there is $q_i\in\matrM_{2M+\ell}(B)$ with $r_i\sim q_i$. Consequently, 
$x_i=[r_i]-M[1_i]$ can be represented as a
difference of projections $[q_i]-M[1_{i}]$ where $q_i\in\matrM_{2M+\ell}(B_i)$.
Defining $q=(q_i)\in\matrM_{2M+\ell}(\myPr)$,
$[q]-M[\unPr]$ is a preimage of $(x_i)$.

If $\eta^0(x)=(x_i)$ and every $x_i$ is positive, we have
$0\leq x_i\leq M[1_i]$
for some fixed $M$. Hence $x_i=[p_i]$ for some $p_i$ which we may
assume lies in $M_{M+\ell}(B_i)$ as above. Consequently $p=(p_i)$
defines an element of $K_0(\myPr)$, and $x=[p]$ by injectivity of
$\eta^0$.
To establish purity when $\valii(\{B_i\})\leq \ell'$,
assume that $x=ny$ in $\prod{K_0(B_i)}$, where 
$x\in\Im\eta^0$ so that for some $M$, $Mn[1_i]\pm x_i\geq 0$.
We conclude that 
${(M+\ell')}[1_{i}]\pm y_i \geq 0$, whence $y\in\Im\eta^0$.
 Proving surjectivity of $\eta^0$ when
$\valiii(\{B_i\})$ is finite is straightforward.
\end{demo}

\begin{lemma}{aboutetaone}
Let  $B_i$ is a sequence of \cstar-algebras.
 If $\rera(\{B_i\})=0$, then
$\eta^1$ is injective. If $\valiv(\{B_i\})<\infty$, then 
$\eta^1$ is surjective.
\end{lemma}
\begin{demo}
To prove injectivity, we assume that $\eta^1(x)=0$ with $x=[(u_i)]$
and $u_i\in \matrM_N(B_i)$. By  \cite{hl:anefscrrz}, $u_i$ is
homotopic to $\mul{n}1_{i}$ within $\UG_n(B_i)$, and since $\matrM_n(B_i)$ has 
finite exponential length because $\rera(B_i)=0$ (see \cite{hl:ercrrzbpc}), we can combine these paths to one from
 $(u_n)$ to the unit of $\UG_n(\myPr)$.
 Proving surjectivity of $\eta^1$ when
$\valiv(\{B_i\})$ is finite is straightforward.
\end{demo}
\begin{lemma}{prosumofadmissible}
Let $B_i$ be a sequence of unital \cstar-algebras and abbreviate
$\Pi=\myPr$, $\Pi/\Sigma=\myPrSu$. Then
\begin{rolist}
\item $\vali(\Pi),\vali(\Pi/\Sigma)\leq\vali(\{B_i\})$.
\item $\valii(\Pi),\valii(\Pi/\Sigma)\leq\valii(\{B_i\})$ if
$\vali(\{B_i\})<\infty$.
\item $\valiii(\Pi),\valiii(\Pi/\Sigma)\leq\valiii(\{B_i\})$ if
$\vali(\{B_i\})<\infty$.
\item $\valiv(\Pi),\valiv(\Pi/\Sigma)\leq\valiv(\{B_i\})$ if
$\rera(\{B_i\})=0$.
\item $\valv(\Pi),\valv(\Pi/\Sigma)\leq\valv(\{B_i\})$ if
$\vali(\{B_i\}),\valii(\{B_i\})<\infty$.
\end{rolist}
\end{lemma}
\begin{demo}
These claims all follow in a straightforward fashion
from the properties of $\eta^*$
established in Lemma \ref{aboutetanil} and \ref{aboutetaone}. We only
prove (v), which is the most involved result. 
Let $x\in
K_0(\myPr)$ and $n$ be given and assume that
$\valv(\{B_i\})\leq\ell$. With $\eta^0(x)=(x_i)$, we can find
$y_i\in K_0(B_i)$ with $-\ell[1_i]\leq y\leq \ell[1_i]$ and
$x_i-y_i\in n K_0(B_i)$. We know from Lemma \ref{aboutetanil} that $(y_i)=\eta^0(y)$ for some $y\in
K_0(\myPr)$ with $-\ell[\unPr]\leq y\leq\ell[\unPr]$. Since
$\eta^0(x-y)\in n\prod K_0(B_i)$ by construction, we further conclude
by purity
that $\eta^0(x-y)\in n\Im\eta^0$ and, since $\eta^0$ is injective,
that $x-y\in n K_0(\myPr)$.

To prove the result for $\Pi/\Sigma$, we note that
\[
\xymatrix{
{0}\ar[r]&{\Ks(\Sigma)}\ar[r]&{\Ks(\Pi)}\ar[r]&{\Ks(\Pi/\Sigma)}\ar[r]&0}
\]
is exact and apply the argument above to a $x\in K_0(\Pi)$ lifting the given
element in $K_0(\Pi/\Sigma)$.
\end{demo} 
\subsubsection*{Algebraically compact $K$-groups}
An abelian group $G$ is \emph{algebraically compact} when $\Pext(-,G)$
vanishes. This class of groups is well studied (cf.\
\cite[VII]{lf:iagI}), and we are going to use the characterization of
it as those groups for which
\begin{rolist}
\item The subgroup $\bigcap_{n\in\NN} nG$ is divisible.
\item $G$ is complete 
\end{rolist}
hold (\cite{ah:acc}). Here completeness refers to the $\ZZ$-adic
topology (cf.\ \cite[7]{lf:iagI}), and completeness does \emph{not}
(as it does in \cite{lf:iagI}!) 
imply any separation properties.

Note that the quantities $\valii$ and $\valv$ make sense for general
ordered abelian groups with order unit. We extend the notions to such
groups and families of them in the obvious way. When $(G_i,1_i)$ is a
family of groups and order units, we set
\[
\prod_b{G_i}=\left.\left\{(g_i)\right|\exists M\in \NN\,\forall i\in
I:-M1_i\leq g_i\leq M1_i\right\}.
\]

\begin{lemma}{bdisalgco}
Whenever $G_i$ is a sequence of abelian groups, then $\prod G_i/\sum
G_i$ is algebraically compact. If, furthermore, all $G_i$ are ordered
with order units, then 
$
\prod_{b}G_i/\sum G_i$ is  algebraically compact provided  that both
$\valii(\{G_i\})$ and $\valv(\{G_i\})$ are finite.
\end{lemma}
\begin{demo}
When $X=\prod G_i/\sum G_i$, then  $X$ is  algebraically compact 
by 
\cite{ah:sfgusrsa}.  When the $G_i$ are also ordered, let
$X_b=\prod_b G_i/\sum G_i$ and consider
$X_b$ as a subgroup of $X$. We are going to prove that  (i) and (ii)
above hold
for $X_b$ from the fact that they hold for $X$.

Since $\valii(\{G_i\})<\infty$, we immediately get that 
if $x\in X_b$ and $x=my$ in $X$, then $y\in X_b$ if $m\not=0$,
by an argument very similar to the one at the end of proof of
 Lemma \ref{aboutetanil}.
Note that this is stronger than the purity of $X_b$  as a subgroup of $X$.
Fix $m\not=0$ and $x\in
\bigcap_{n\in\NN} nX_b$ and write $x=my$ for $y\in\bigcap_{n\in\NN}
nX$. Applying this observation twice, we get 
\[
y\in X_b\cap \bigcap_{n\in\NN} nX=\bigcap_{n\in\NN}(nX\cap X_b)=\bigcap_{n\in\NN}nX_b
.\]

It remains to show that $X_b$ is complete.
To do so, we first note that for any given $x\in X$ there is $y\in
X_b$ such that $x-y\in\bigcap_{n\in \NN}nX$. For when
$\valv(\{G_i\})\leq\ell$ and  $x=(x_i)+\sum G_i$ is given,
we may 
find $y_i$ such that $-\ell 1_i\leq y_i\leq \ell 1_i$ and $x_i-y_i\in
i! G_i$. Clearly $y=(y_i)+\sum G_i$ has the desired property.
Now let $(x_n)$ be Cauchy in $X_b$, and recall that it converges to some
$x \in X$. With $y$ chosen as above, $(x_n)$ also converges to $y\in X_b$
in the  $\ZZ$-adic topology of $X_b$.
This is because  the $\ZZ$-adic topology of $X_b$ coincides with the topology of
$X_b$ induced by  the $\ZZ$-adic topology of $X$, since $X_b$ is a pure subgroup of $X$.
\end{demo}
\begin{corol}{isac}
Let $B_i$ be a sequence of unital \cstar-algebras.
\begin{rolist}
\item 
$K_1(\myPrSu)$ is algebraically compact if  $\rera(\{B_i\})=0$ and $\valiv(\{B_i\})<\infty$.
\item $K_0(\myPrSu)$ is algebraically compact if either
$\vali(\{B_i\}),\valiii(\{B_i\})<\infty$
or $\vali(\{B_i\}),\valii(\{B_i\}),\valv(\{B_i\})<\infty$. 
\end{rolist}
\end{corol}
\begin{demo}
Part (i) follows from Lemmas \ref{aboutetaone} and \ref{bdisalgco}.
Part (ii) follows from Lemmas \ref{aboutetanil} and \ref{bdisalgco}.
\end{demo}

\subsubsection*{Admissible targets}
We are now ready to collect our results in the special case of
admissible target algebras.
\begin{theor}{collect}
Let $B_i$ be a sequence of admissible target algebras of the same type and let $A$ be
any \cstar-algebra. Then
\begin{rolist}
\item $\ETA:\hugePK{\mysmPr}\longrightarrow \prod{\hugePK{B_i}}$ is
injective.
\item The natural map $\Hom_\Lambda(\hugePK{A},\hugePK{\mysmPr})
\longrightarrow \prod{\Hom_\Lambda(\hugePK{A},\hugePK{B_i})}$ is injective.
\item The natural map $\KK(A,\mysmPrSu)\longrightarrow
\Hom_\Lambda(\hugePK{A},\hugePK{\mysmPrSu})$ is an isomorphism if $A$ satisfies the UCT.
\item $\myPr$ and $\myPrSu$ are admissible targets
\end{rolist}
\end{theor}
\begin{demo}
For (i), decompose $\ETA$ into maps $\etaKs$ and $\eta^*_n$. We get from Lemma \ref{aboutetanil} and \ref{aboutetaone} that
$\eta^*$ are injections whose images are pure subgroups. Injectivity
of $\eta^*_n$ then follows
by a diagram chase on
\[
\xymatrix{
{K_{*}(\myPr)}\ar[r]^-{\times n}\ar[d]_{\eta^{*}}&
{K_{*}(\myPr)}\ar[r]^-{\rho^{*}_n}\ar[d]_{\eta^{*}}&
{\pKstar{n}{\myPr}}\ar[r]^-{\beta^{*}_n}\ar[d]_{\eta^{*}_n}&
{K_{*+1}(\myPr)}\ar[d]_{\eta^{*+1}}\\
{\prod{K_{*}(B_i)}}\ar[r]_-{\times n}&
{\prod{K_{*}(B_i)}}\ar[r]_-{\prod\rho^{*}_n}&
{\prod{\pKstar{n}{B_i}}}\ar[r]_-{\prod\beta^{*}_n}&
{\prod{K_{*+1}(B_i).}}
}
\]
Claim (ii) is a direct
consequence of (i), and (iii) follows by combining the UMCT \ref{UMCT}
 with Corollary
\ref{isac}. Finally, (iv) follows by Lemma \ref{prosumofadmissible}.
\end{demo}

\subsection{Partial maps on \underline{\textbf{K}}$({-})$}\label{partialmaps}
In this appendix we concern ourselves with associating $K$-theoretical
data to completely positive contractions. Starting from such maps, say
from $A$ to $B$, we
shall be able to induce maps sending a finite set of projections
representing a finite part of the $K$-theory of $A$ to elements of
the $K$-theory of $B$. 

Although there are advantages of doing this even for subsets
representing elements of $K_0(A)$, the real strength of this approach
only surfaces when we work with all of $\hugePK{A}$ and $\hugePK{B}$. Our \emph{partial maps} do not descend to
well-defined maps on subsets of $\hugePK{A}$, let alone to all of
$\hugePK{A}$, but this fact
 does not cause any problems except
notational and technical inconveniences. 

As noted in \cite{mdgg:crahcrrz} we can realize any element of $\hugePK{-}$ as a difference of
projections from
\[
\underline{\bigproj}(A)=\bigcup_{m\geq 1} \bigproj(A\otimes C({\mathbb T})\otimes
C(W_m)\otimes\KKK)
\]
where the $W_m$ are the Moore spaces of order $m$. This picture of
$\hugePK{A}$ encompasses the standard pictures of $K_0(A)$ and
$\Ks(A)$ using projections of $A\otimes\KKK$ and $A\otimes
C({\mathbb T})\otimes\KKK$, respectively, but not the standard picture of
$K_1(A)$ using unitaries of $(A\otimes\KKK)^{\sim}$. We need to pay
special attention to this.  Checking the facts stated as lemmas
below is tedious but straight-forward. We leave it to the reader with
due apologies.

\begin{defin}{Ktriples}
Let $A$ be a \cstar-algebra. A \emph{$\huK$-triple
$(\projsetKhu,\fseta,\delta)$} consists of finite subsets
$\projsetKhu\subseteq \underline{\bigproj}(A)$ and $\fseta\subseteq A$ and a
$\delta>0$ chosen such that whenever $\phi$ is a completely positive
contraction which is $\delta$-multiplicative on $\fseta$, then
\[
\tfrac12 \not\in\spectrum((\phi\otimes \id)(p))
\]
for each $p\in\projsetKhu$, where $\id$ is the identity of $C(\mathbb{T})\otimes
C(W_m)\otimes\KKK$ for suitable $m$. A \emph{$K_1$-triple}
$(\uniset,\fseta,\delta)$ consists of finite subsets
$\uniset\subseteq \UG((A\otimes\KKK)^\sim)$ and $\fseta\subseteq A$ and a
$\delta>0$ chosen such that whenever $\phi$ is a completely positive
contraction which is $\delta$-multiplicative on $\fseta$, then
\[
0\not\in\spectrum((\phi\otimes \id)(v))
\]
for each $v\in\uniset$.
\end{defin}

We define $K_0$- and $\Ks$-triples analogously
to the $\huK$-triple case, by using projections in $A \otimes \KKK$ and
$A \otimes C(\mathbb T) \otimes \KKK$, respectively.

The following lemma shows that any finite subset of projections or
unitaries can be \emph{complemented} to a triple of the appropriate kind.

\begin{lemma}{getdefined}
Let $A$ and $C$ be \cstar-algebras, and fix finite sets
$\projset,\uniset\subseteq (A\otimes C)^{\sim}$ consisting of projections and
unitaries, respectively. Then there
exists $\delta>0$ and a finite set $\fseta\subseteq A$ such that
whenever $\phi:A\longrightarrow B$ is a unital completely positive map
which is $\delta$-multiplicative on $\fseta$, then
\[
\tfrac12\not\in\spectrum((\phi\otimes \id_C)(p))\qquad
0\not\in\spectrum((\phi\otimes \id_C)(v))
\]
for every $p\in\projset$, $v\in\uniset$.
\end{lemma}

Let $\chi_0:[0,\tfrac12)\cup(\tfrac12,1]\rightarrow [0,1]$ be $0$ on
$[0,\tfrac 12)$ and $1$ on $(\tfrac12,1]$, and
let $\chi_1:{(0,1]}\rightarrow {\mathbb C}$ be
given by $\chi_1(x)=x^{-1/2}$.

\begin{defin}{thepartialmaps}
Let $(\projset,\fseta,\delta)$ be a $\huK$-triple, and assume that
$\phi:A\longrightarrow B$ is a completely positive contraction which is
$\delta$-multiplicative on $\fseta$. We define
$\phi\pstar:\projset\longrightarrow \hugePK{B}$ by
\[
\phi\pstar(p)=[\chi_0(\phi\otimes\id)(p)].
\]
When $(\uniset,\fseta,\delta)$ is a
$\huK$-triple, we define
$\phi\pstar:\uniset\rightarrow K_1(B)$ by
\[
\phi\pstar(v)=[V\chi_1(VV^*))]
\]
where $V=(\phi \otimes \id)(v)$.
\end{defin}
Maps into $K_0(B)$ and $\Ks(B)$ are defined from $K_0$-triples and
$\Ks$-triples similarly to the $\huK$-triple case.

\begin{lemma}{aha}
Let $A$ be a unital \cstar-algebra and $(\projset,\fseta,\delta)$ a
$\huK$-triple. Let $\phi:A\longrightarrow B$ be a completely positive
map which is $\delta$-multiplicative on $\fseta$ and let
$j:B\longrightarrow C$ be a unital $*$-homomorphism. Then
$
(j\phi)\pstar(p)=j_*\phi\pstar(p)
$
for all $p\in\projset$.
\end{lemma}

To establish the next two results, one may use that the
canonical map from $K_0(A)$ to $K_1(SA)$ is defined using \emph{scalar}
rotation matrices. 
\begin{lemma}{star01}
Whenever a $\Ks$-triple $(\projsetKs,\fseta,\delta)$ is given, there
exist a $K_0$-triple $(\projsetnil,\fseta_0,\delta_0)$ and a
$K_1$-triple $(\uniset,\fseta_1,\delta_1)$ with
$\delta_i<\delta$ and $\fseta_i\supseteq \fseta$ such that if
for two unital completely positive contractions $\phi,\psi$ which are
$\delta_i$-multiplicative on $\fseta_i$,
\[
\phi\pstar(p)=\psi\pstar(p)\qquad
\phi\pstar(v)=\psi\pstar(v)
\]
for all $p\in\projsetnil$ and all $v\in\uniset$, then
$
\phi\pstar(p)=\psi\pstar(p)
$
for all $p\in\projsetKs$.
\end{lemma}
\begin{lemma}{01star}
Whenever a $K_0$-triple $(\projsetnil,\fseta_0,\delta_0)$ and a
$K_1$-triple $(\uniset,\fseta_1,\delta_1)$ is given, there exists a
$\Ks$-triple $(\projsetKs,\fseta,\delta)$ with
$\delta<\delta_i$ and $\fseta\supseteq \fseta_i$ such that if
for two unital completely positive contractions $\phi,\psi$ which are
$\delta$-multiplicative on $\fseta$,
$
\phi\pstar(p)=\psi\pstar(p)
$
for all $p\in\projsetKs$, then
\[
\phi\pstar(p)=\psi\pstar(p)\qquad
\phi\pstar(v)=\psi\pstar(v)
\]
for all $p\in\projsetnil$ and all $v\in\uniset$.
\end{lemma}

\begin{lemma}{close}
Let $(\projset,\fseta,\delta)$ be a $\huK$-triple on $A$. There exists
$\epsilon>0$ and a finite set $\fset\subseteq A$ such that if $\phi,\psi$ are
completely positive contractions which are $\delta$-multiplicative on
$\fseta$ and satisfy $\norm{\phi(a)-\psi(a)}<\epsilon$ for all $a\in\fset$, then 
$
\phi\pstar(p)=\psi\pstar(p)
$
for all $p\in\projset$.
\end{lemma}


\begin{thebibliography}{Kas80b}

\bibitem[BK97]{bbek:gilfdc}
B.~Blackadar and E.~Kirchberg, \emph{Generalized inductive limits of
  finite-dimensional ${C}^*$-algebras}, Math. Ann. \textbf{307} (1997), no.~3,
  343--380.

\bibitem[Bla86]{bb:koa}
B.~Blackadar, \emph{{$K$}-theory for operator algebras}, Math. Sci. Research
  Inst. Publ., vol.~5, Springer-Verlag, New York, 1986.
 
\bibitem[Bra72]{ob:af}
O.~Bratteli, \emph{Inductive limits of finite-dimensional \cstar-algebras},
Trans. Amer. Math. Soc. \textbf{171} (1972), 195--234.

\bibitem[BDF77]{bdf:eck}
L.G. Brown, R.G. Douglas and P.A. Fillmore, \emph{Extensions of
{$C^*$}-algebras and {$K$-homology}}, Ann. of Math. \textbf{105} (1977), no.~2, 265--324.

\bibitem[BP91]{lgbgkp:crrz}
L.G. Brown and G.K. Pedersen, \emph{{$C^*$}-algebras of real rank zero}, J.\
  Funct.\ Anal. \textbf{99} (1991), 131--149.

\bibitem[Bro77]{lgb:sihsc}
L.G. Brown, \emph{Stable isomorphism of hereditary subalgebras of
  {$C^*$}-algebras}, Pacific J.\ Math. \textbf{71} (1977), no.~2,
  335--348.


\bibitem[CE76]{mdcege:cplpc}
M.D. Choi and E.G. Effros, \emph{The completely positive lifting problem for
  {$C^*$}-algebras}, Ann. of Math. (2) \textbf{104} (1976), no.~3, 585--609.

\bibitem[CH90]{acnh:dmakb}
A.~Connes and N.~Higson, \emph{D\'eformations, morphismes asymptotiques et
  ${K}$-th\'eorie bivariante}, C. R. Acad. Sci. Paris S\'er. I Math.
  \textbf{311} (1990), no.~2, 101--106.

\bibitem[Cun81]{jc:kcc}
J.~Cuntz, \emph{{$K$}-theory for certain {$C^*$}-algebras}, Ann.\ of Math.
  \textbf{113} (1981), 181--197.

\bibitem[Dar97]{md:nsoaf}
M.~{D\u{a}d\u{a}rlat}, \emph{Nonnuclear subalgebras of AF algebras},
preprint 1997.

\bibitem[DG97]{mdgg:crahcrrz}
M.~{D\u{a}d\u{a}rlat} and G.~Gong, \emph{A classification result for
  approximately homogeneous ${C}^*$-algebras of real rank zero}, Geom. Funct.
  Anal. \textbf{7} (1997), no.~4, 646--711.

\bibitem[DL92]{mdtal:kasa}
M.~{D\u{a}d\u{a}rlat} and T.A. Loring, \emph{The {$K$}-theory of {Abelian}
  subalgebras of {$AF$} algebras}, J.\ reine angew.\ Math. \textbf{432} (1992),
  39--55.

\bibitem[DL94]{mdtal:ecrrzc}
M.~{D\u{a}d\u{a}rlat} and T.A. Loring, \emph{Extensions of certain real rank
  zero ${C}^*$-algebras}, Ann.\ Inst.\ Fourier \textbf{44} (1994), 907--925.

\bibitem[DL96a]{mdtal:ccomk}
M.~D\u{a}d\u{a}rlat and T.A. Loring, \emph{Classifying {$C^*$}-algebras via
  ordered, mod-$p$ ${K}$-theory}, Math. Ann. \textbf{305} (1996), no.~4,
  601--616.

\bibitem[DL96b]{mdtal:umctkg}
M.~{D\u{a}d\u{a}rlat} and T.A. Loring, \emph{A universal multicoefficient
  theorem for the {K}asparov groups}, Duke Math. J. \textbf{84} (1996), no.~2,
  355--377.

\bibitem[EG96]{gaegg:ccrrzii}
G.A. Elliott and G.~Gong, \emph{On the classification of {$C^*$}-algebras of
  real rank zero, {II}}, Ann.\ of Math. \textbf{144} (1996), 497--610.

\bibitem[EL]{setal:ccsr}
S.~Eilers and T.A. Loring, \emph{Computing contingencies for stable relations},
  \emph{Internat.\ J.\ Math.}, to appear.

\bibitem[Ell76]{gae:cilssfa}
G.A. Elliott, \emph{On the classification of inductive limits of sequences of semisimple
finite-dimensional algebras}, J.\ Algebra \textbf{38} (1976), 29--44.

\bibitem[Ell93]{gae:ccrrz}
G.A. Elliott, \emph{On the classification of {$C^*$}-algebras of real rank
  zero}, J.\ reine angew.\ Math. \textbf{443} (1993), 179--219.

\bibitem[Fu70]{lf:iagI}
L.\ Fuchs, \emph{Infinite abelian groups {I}}, Academic Press, New York, San
  Francisco, London, 1970.

\bibitem[Hig87a]{nh:ck}
N.~Higson, \emph{A characterization of {$\KK$}-theory}, Pacific J. Math.
  \textbf{126} (1987), no.~2, 253--276.

\bibitem[Hig87b]{nh:tt}
N.~Higson, \emph{On a technical theorem of Kasparov  }, 
J. Funct. Anal. 
  \textbf{73} (1987), no.~1, 107--112

\bibitem[Hig90]{nh:cfek}
N.~Higson, \emph{Categories of fractions and excision in {$\KK$}-theory}, J.
  Pure Appl. Algebra \textbf{65} (1990), no.~2, 119--138.

\bibitem[Hig95]{nh:cetd}
N.~Higson, \emph{{$C^*$}-algebra extension theory and duality}, J. Funct. Anal.
  \textbf{129} (1995), no.~2, 349--363.

\bibitem[HLT96]{thkt:uct}
T.G. Houghton-Larsen and K.~Thomsen, \emph{Universal (co)homology theories},
  preprint, 1996.

\bibitem[Hul62a]{ah:acc}
A.~Hulanicki, \emph{On algebraically compact groups}, Bull. Acad. Polon. Sci.
  S\'er. Sci. Math. Astronom. Phys. \textbf{10} (1962), 71--75.

\bibitem[Hul62b]{ah:sfgusrsa}
A.~Hulanicki, \emph{The structure of the factor group of the unrestricted sum
  by the restricted sum of {A}belian groups}, Bull. Acad. Polon. Sci. S\'er.
  Sci. Math. Astronom. Phys. \textbf{10} (1962), 77--80.

\bibitem[JT91]{kkjkt:ek}
K.K. Jensen and K.~Thomsen, \emph{Elements of {$K\!K$}-theory},
  {Birkh\"{a}user}, Boston, 1991.

\bibitem[Kas80a]{ggk:hcmtst}
G.G. Kasparov, \emph{Hilbert {$C^*$}-modules: {Theorems} of {S}tinespring and
  {V}oiculescu}, J. Operator Theory \textbf{4} (1980), no.~1, 133--150.

\bibitem[Kas80b]{ggk:okec}
G.G. Kasparov, \emph{The operator ${K}$-functor and extensions of
  {$C^*$}-algebras}, Izv. Akad. Nauk SSSR Ser. Mat. \textbf{44} (1980), no.~3,
  571--636, 719.

\bibitem[Kir94]{ek:cpicukt}
E.\ Kirchberg, \emph{The classification of purely infinite {$C^*$}-algebras
  using {Kasparov's} theory}, preprint, third draft, 1994.

\bibitem[KP98]{ekncp:eeccfco}
E.\ Kirchberg and N.C. Phillips, \emph{Embedding of exact {$C^*$}-algebras and
  continuous fields in the {C}untz algebra ${{\mathcal O}}_2$}, preprint, 1998.

\bibitem[Lin93]{hl:ercrrzbpc}
H.~Lin, \emph{Exponential rank of {$C^*$}-algebras with real rank zero and the
  {Brown-Pedersen} conjectures}, J.\ Funct.\ Anal. \textbf{114} (1993), no.~1,
  1--11.

\bibitem[Lin96]{hl:anefscrrz}
H.\ Lin, \emph{Approximation by normal elements with finite spectra in
  ${C}^*$-algebras of real rank zero}, Pacific J. Math. \textbf{173} (1996),
  no.~2, 443--489.

\bibitem[Lin97]{hl:saueh}
H.\ Lin, \emph{Stably approximately unitary equivalence of homomorphisms},
  preprint, 1997.

\bibitem[Lin98a]{hl:newtaf}
H.\ Lin, \emph{Classification of simple tracially {AF} {$C^*$}-algebras},
  preprint, 1998.

\bibitem[Lin98b]{hl:taf}
H.\ Lin, \emph{Tracially {AF} {$C^*$}-algebras}, preprint, 1998.

\bibitem[Lor97]{tal:lsppc}
T.A. Loring, \emph{Lifting solutions to perturbing problems in
  {$C^*$}-algebras}, Fields Institute Monographs, vol.~8, American Mathematical
  Society, Providence, RI, 1997.

\bibitem[LP95]{hlncp:ammco}
H.X. Lin and N.C. Phillips, \emph{Almost multiplicative morphisms and the
  {C}untz algebra $ {{\mathcal O}}_2$}, Internat. J. Math. \textbf{6} (1995),
  no.~4, 625--643.

\bibitem[Pas81]{wlp:kcca}
W.L. Paschke, \emph{${K}$-theory for commutants in the {C}alkin algebra},
  Pacific J. Math. \textbf{95} (1981), no.~2, 427--434.

\bibitem[Ped79]{gkp:cag}
G.K. Pedersen, \emph{{$C^*$}-algebras and their automorphism groups}, Academic
  Press, London, 1979.

\bibitem[Phi94]{ncp:ctnpisc}
N.C. Phillips, \emph{A classification theorem for nuclear purely infinite
  simple {$C^*$}-algebras}, Preprint, 1994.

\bibitem[Phi97]{ncp:auehoca}
N.C. Phillips, \emph{Approximate unitary equivalence of homomorphisms from odd
  {C}untz algebras}, Operator algebras and their applications (Waterloo, ON,
  1994/1995), Fields Inst. Commun., vol.~13, Amer. Math. Soc., Providence, RI,
  1997, pp.~243--255.

\bibitem[Pop97]{sp:lfdac}
S.~Popa, \emph{On local finite-dimensional approximation of ${C}\sp
  *$-algebras}, Pacific J. Math. \textbf{181} (1997), no.~1, 141--158.

\bibitem[Rie83]{mr:dsrkc}
M.A. Rieffel, \emph{Dimension and stable rank in the {$K$}-theory of
  {$C^*$}-algebras}, Proc.\ London Math.\ Soc. \textbf{46} (1983), no.~3,
  301--333.

\bibitem[Rie87]{mar:hgugnt}
M.A. Rieffel, \emph{The homotopy groups of the unitary groups of noncommutative
  tori}, J. Operator Theory \textbf{17} (1987), no.~2, 237--254.

\bibitem[R{\o}r93]{mr:cilca}
M.~R{\o}rdam, \emph{Classification of inductive limits of {Cuntz} algebras},
  J.\ reine angew.\ Math. \textbf{440} (1993), 175--200.

\bibitem[R{\o}r95]{mr:ccisa}
M.~R{\o}rdam,  \emph{Classification of certain infinite simple 
\cstar-algebras}, J. Funct. Anal. \textbf{131}
(1995), no.~2, 415--458


\bibitem[RS87]{jrcs:ktuctkgk}
J.~Rosenberg and C.~Schochet, \emph{The {K\"u}nneth theorem and the universal
  coefficient theorem for {K}asparov's generalized {$K$}-functor}, Duke Math.\
  J. \textbf{55} (1987), 431--474.



\bibitem[Sal92]{ns:rqaKKt}
N.\ Salinas,  \emph{Relative quasidiagonality and $KK$-theory}, Houston J. Math. \textbf{18} (1992), no.~1, 97--116

\bibitem[Sch84]{cs:tmcIV}
C.\ Schochet, \emph{Topological methods for {$C^*$}-algebras {IV}: Mod $p$
  homology}, Pacific J.\ Math. \textbf{114} (1984), 447--468.

\bibitem[Sch84]{cs:fskg}
C.\ Schochet, \emph{The fine structure of the Kasparov groups II: relative
quasidiagonality}, Preprint.

\bibitem[Ska88]{gs:nnk}
G.~Skandalis, \emph{Une notion de nucl\'earit\'e en ${K}$-th\'eorie (d'apr\`es
  {J}.\ {C}untz)}, $K$-Theory \textbf{1} (1988), no.~6, 549--573.

\bibitem[Val83]{av:rkg}
A.~Valette, \emph{A remark on the {K}asparov groups
  {$\operatorname{Ext}({A},\,{B})$}}, Pacific J. Math. \textbf{109} (1983),
  no.~1, 247--255.

\bibitem[Voi76]{dv:nwvt}
D.~Voiculescu, \emph{A non-commutative {W}eyl--von {N}eumann theorem}, Rev.
  Roumaine Math. Pures Appl. \textbf{21} (1976), no.~1, 97--113.

\bibitem[Voi93]{dv:aqo}
D.~Voiculescu, \emph{Around quasidiagonal operators}, Integral Equations
  Operator Theory \textbf{17} (1993), no.~1, 137--149.

\bibitem[Zha90]{sz:ppisc}
S.~Zhang, \emph{A property of purely infinite simple ${C}^*$-algebras}, Proc.
  Amer. Math. Soc. \textbf{109} (1990), no.~3, 717--720.

\end{thebibliography}
\providecommand{\bysame}{\leavevmode\hbox to3em{\hrulefill}\thinspace}

\begin{center}
\parbox{6.5cm}
{\begin{center}{\sc
Marius D\u{a}d\u{a}rlat\\
Department of Mathematics\\
Purdue University\\
West Lafayette\\
IN 47907\\
U.S.A.\\}
{\tt mdd\atsign math.purdue.edu}
\end{center}}
\makebox[0.5cm]{}
\parbox{6.5cm}
{\begin{center}{\sc
S\o ren Eilers\\
Matematisk Afdeling\\
K\o benhavns Universitet\\
Universitetsparken 5\\
DK-2100 Copenhagen \O\\
Denmark\\}
{\tt eilers\atsign math.ku.dk}
\end{center}}\\
\end{center}
\end{document}